\documentclass[twoside,11pt]{article}
\usepackage{amsmath, amsthm, amsfonts, amssymb}
\usepackage{mathrsfs}
\def\Bbb{\mathbb}
\newtheorem{theorem}{Theorem}[section]
\newtheorem{lemma}[theorem]{Lemma}
\newtheorem{corollary}[theorem]{Corollary}
\newtheorem{definition}[theorem]{Definition}
\newtheorem{remark}[theorem]{Remark}

\makeatletter \@addtoreset{equation}{section} \makeatother

\def\dmin{\displaystyle\min}
\def\dbigvee{\displaystyle\bigvee}
\def\dsum{\displaystyle\sum}
\def\dint{\displaystyle\int}
\def\dmax{\displaystyle\max}

\def\dsup{\displaystyle\sup}
\def\dinf{\displaystyle\inf}

\begin{document}
\baselineskip=13pt
\renewcommand{\arraystretch}{2}
\arraycolsep=1.6pt

\title{\Large \bf Comparison
Theorems of Kolmogorov Type  for Classes\\ Defined by Cyclic
Variation Diminishing Operators\\ and Their Application
\thanks
{Project   supported by the National Natural Science Foundation of
China(Grant No. 10371009) and Research Fund for the Doctoral
Program Higher Education(Grant No. 20050027007).
\newline{$^{**}$ Corresponding author.}
\newline {\small \indent \ E-mail address: fanggs@bnu.edu.cn}}}
\author{ $^{**}$Fang Gensun\ \quad Li Xuehua\vspace{-0.3cm}\\
\vspace{-0.4cm} {\small School of Mathematical Sciences,
Beijing Normal University,}
{\small  Beijing 100875, China}\\
\\}
\date{ }
\maketitle
\vspace{-18mm}

{\begin{abstract}Using present a unified approach, we establish a
Kolmogorov type comparison theorem  for the  classes of
$2\pi$-periodic functions defined by a special class of operators
having certain oscillation properties, which includes the
classical Sobolev class of functions with 2$\pi$-periodic, the Achieser class,
and the Hardy-Sobolev class as its special examples.  Then, using these
results, we prove a Taikov type inequality, and calculate the
exact values of Kolmogorov, Gel$'$fand, linear and information
$n$--widths of this class of functions in some  space $L_{q}$,
which is the classical Lebesgue integral space of 2$\pi$--periodic
with the usual norm.

\vspace{0.2cm} { \bf Key words and phrases:} Comparison theorems of
Kolmogorov type,  $n$-widths, oscillation properties

\vspace{0.2cm} {\bf AMS classification:} 41A46, 30D55
\end{abstract}

\section{Introduction}\label{sec1}

\hspace{1ex} Osipenko \cite{Oc97b} noticed that many extremal
problems of approximation theory, such as optimal recovery,
optimal quadrature formulae, and  n-widths
etc., can be solved for the classes of smooth periodic functions
which can be represented as convolutions with cyclic variation
diminishing kernels. However, the functions in some classes of
analytic functions do not seem to be representable as convolutions
with cyclic variation diminishing kernels. Nevertheless, results
that have been known for these classes are very close to those
obtained in the smooth case above (see
\cite{Oc94,Oc95b,Oc97a,FGSLXH,FGSLXH2}). The idea of constructing
general theory, which is covering both of the smooth and analytic
cases, has been repeatedly expressed by Tikhomirov (see
\cite{Oc97b}). Osipenko \cite{Oc97b} proposed a unified approach
to the problem of the exact calculation of $n$-widths on the
classes of functions, defined by a  special class of (in general,
nonlinear) operators that have the property of cyclic variation
diminishing. In this paper, we continue the work of Osipenko
\cite{Oc97b} to study comparison theorems of Kolmogorov type,
inequalities of Taikov type, and determine the precise values $n$--widths
on such class of functions.

 First we recall some definitions.
Let $W$ be a subset of a linear space $X$. We
 consider the problem of the optimal recovery of a linear
 functional $L$ on this subset from the values of linear
 functionals $l_{1},\ldots,l_{n}$. For $x\in W$, we set
 $$
 Ix:=(l_{1}x,\ldots,l_{n}x).
 $$
 The operator $I$:$W\rightarrow K^{n}$, where $K={\mathbb R}$ or
 ${\mathbb C}$, depending on whether $X$ is a real or a complex
 space, is called the information operator. The quantity
 $$
 e(L,W,I):=\inf_{S: K^{n}\rightarrow K}\sup_{x\in W}|Lx-S(Ix)|
 $$
 is called the error of the optimal recovery of the functional $L$
 on the set $W$. Each method $S_{0}$ such that
 $$
 \sup_{x\in W}|Lx-S_{0}(Ix)|=e(L,W,I)
 $$
 is said to be an optimal method of recovery.
 Smolyak \cite{Sm} proved that in the real case, for a
 convex centrally symmetric set $W$, there exists among optimal
 methods of recovery a linear method and
 \begin{equation}
 e(L,W,I):=\sup_{x\in W,\,Ix=0}|Lx|.
  \label{dt6-1}
 \end{equation}
Each element $x_{0}\in W$ such that $Ix_{0}=0$ and
 $$
 |Lx_{0}|=\sup_{x\in W,\,Ix=0}|Lx|,
 $$
 is said to be extremal. The problem of finding an extremal element
 often turns out to be simpler than that of finding an optimal recovery method. \\

 Let $X$ be a normed linear space
with norm $\|\cdot\|$, and $X_{n}$ be an
$n$--dimensional subspace of $X$. For each $x\in X$, $E(x;X_{n})$
denotes the distance of the $n$--dimensional subspace $X_{n}$ from
$x$, defined by
\begin{equation}
 E(x;X_{n}):=\inf_{y\in X_n}\|x-y\|,
\label{dt4-1}
\end{equation}
and the quantity
\begin{equation}
 E(A,X_{n}):=\sup_{x\in A}\inf_{y\in X_n}\|x-y\|
\label{dt4-2}
\end{equation}
is said to be the deviation of $A$ from $X_n$. Thus $E(A,X_{n})$
measures how well the ``worst element" of $A$ can be approximated
from $X_{n}$.\\

 Given a subset $A$ of $X$, one might ask how well one
can approximate $A$ by $n$--dimensional subspaces of $X$. Thus, we
consider the possibility of allowing the $n$--dimensional
subspaces $X_{n}$ to vary within $X$.  This idea, introduced by
Kolmogorov in 1936,  is now referred to as the Kolmogorov
$n$--width of $A$ in $X$. It is defined by
\begin{equation}
 d_n(A,X):=\inf_{X_n}\sup_{x\in A}\inf_{y\in X_n}\| x-y\|,
\label{dt4-3}
\end{equation}
where $X_n$ runs over all $n$--dimensional subspaces of $X$.

 The Kolmogorov $n$--width $d_n(A,X)$ describes the
minimum error of $A$ approximated by any $n$--dimensional subspace
$X_n$ in $X$. In addition to the Kolmogorov $n$--width, there are
three other related concepts that will be studied in this paper.
The linear $n$--width of $A$ in $X$ is defined by
\begin{equation}
\lambda_n(A,X):=\inf_{P_n}\sup_{x\in A}\| x-P_{n}x\|,
\label{dt4-4}
\end{equation}
where the infimum is taken over all  bounded linear operators
mapping $X$ into itself whose range has dimension at most $n$.

 The Gel$'$fand $n$--width of $A$ in $X$ is given by
\begin{equation}
d^n(A,X):=\inf_{X^n}\sup_{x\in A\cap X^n}\| x\|, \label{dt4-5}
\end{equation}
where $X^{n}$ runs over all subspaces of $X$ of codimension $n$
(here we assume that $0\in A$), and  the information n--width is
the quantity
\begin{equation}
i_{n}(A,X):=\inf_{l_{1},\ldots,l_{n}}\inf_{m:Z^{n}\rightarrow X}
\sup_{x\in A}\|x-m(l_{1}x,\ldots,l_{n}x)\|, \label{dt4-6}
\end{equation}
where $l_{1},\ldots,l_{n}$ run over all continuous linear
functionals on $A$ and $m$ is taken over all maps of $Z^{n}$ into
$X$ ($Z={\mathbb R}$ or ${\mathbb C}$ depending on whether $A$ is
a set of real--valued or complex--valued functions). Fore more detailed
information about the $n$--widths, we refer to monographs  of
Pinkus \cite{P}, and Lorentz, Golischek and Makovoz \cite{Lo}.

Now we recall some notions of sign changes of vectors
and functions, and  the definitions of Property $B$  and $NCVD$.

\begin{definition}\label{Definition 1}
(see \cite{P}, pp. 45, 59) Let $x=(x_1,\ldots,x_n)\in{\mathbb
R}^n\backslash\{0\}$ be a real non--trivial
vectors. \\
\noindent (i) $S^{-}(x)$ indicates the number of sign changes in the
sequence $x_1,\ldots,x_n$ with zero terms discarded. The number
$S_{c}^{-}(x)$ of cyclic variations of sign of $x$ is given by
$$
S_{c}^{-}(x):=\max_{i}
S^{-}(x_i,x_{i+1},\ldots,x_n,x_1,\ldots,x_i)
=S^{-}(x_k,\ldots,x_n,x_1,\ldots,x_k),
$$
where $k$ is any integer for which $x_k\neq 0$. Obviously
$S_{c}^{-}(x)$ is invariant
under cyclic permutations, and $S_{c}^{-}(x)$ is always an even number.
\\

\noindent (ii) $S^{+}(x)$ counts the maximum number of sign changes in the
sequence $x_1,\ldots,x_n$ where zero terms are arbitrarily
assigned values $+1$ or $-1$. The number $S^+_{c}(x)$ of maximum
cyclic variations of sign of $x$ is defined by
$$
S_{c}^{+}(x):=\max_{i}
S^{+}(x_i,x_{i+1},\ldots,x_n,x_1,\ldots,x_i).
$$

 Let $f$ be a piecewise continuous, $2\pi$--periodic
function. We assume that $f(x)=[f(x+)+f(x-)]/2$ for all $x$ and
\begin{equation}
\ \  S_{c}(f):=\sup
S_{c}^{-}\left((f(x_1),\ldots,f(x_m))\right),\label{s2}
\end{equation}
where the supremum is taken over all $x_1<\cdots<x_m<x_1+2\pi$ and
all $m\in{\mathbb N}$.
\end{definition}

 Suppose that $f$ is a continuous function of periodic
$2\pi$. We define
\begin{equation} \widetilde{Z}_{c}(f):=\sup
S_{c}^{+}\left((f(x_1),\ldots,f(x_m))\right), \label{s3}
\end{equation}
where the supremum runs over all $x_1<\ldots<x_m<x_1+2\pi$ and all
$m\in {\mathbb N}$.   Clearly,
$$
S_{c}(f)\leq\widetilde{Z}_c(f),
$$
where $S_{c}(f)$ denotes the number of sigh changes of $f$ on a period,
$\widetilde{Z}_c(f)$ denotes the number of zeros of $f$ on a
period, and  sign changes are counted once and zeros which are
not sign changes are counted twice.

\begin{definition}\label{Definition 2}(see \cite{P}, p.129) A
real, $2\pi$--periodic, continuous function $G$ satisfies Property
$B$ if for every choice of $0\leq t_{1}<\cdots<t_{m}<2\pi$ and
each $m\in {\mathbb N}$, the subspace
$$
X_{m}:=\{b+\sum_{j=1}^{m}b_{j}G(\cdot-t_{j}):\sum_{j=1}^{m}b_{j}=0\}
$$
is of dimension $m$, and is a weak Tchebycheff ($WT$--) system
(see \cite{P}, p.39) for all $m$ odd.
\end{definition}

 Let $\phi$ be a piecewise continuous $2\pi$--periodic
function satisfying $\phi\perp{\bf 1}$ and set
$\psi(x):=a+(G\ast\phi)(x)$. If $G$ satisfies Property $B$, then
$S_c(\psi)\leq S_c(\phi)$(see \cite[p.129]{P}).

\begin{definition}\label{Definition 3} (see \cite{P}, p.60 and
p.126)\ \  We say that a real, continuous, $2\pi$--periodic
function $k$ is non--degenerate cyclic variation diminishing
($NCVD$) if $S_c(k\ast h)\leq S_c(h)$ for each real,
$2\pi$--periodic, piecewise continuous function $h$, and
$$
\dim{\rm span}\{k(x_{1}-\cdot),\ldots,k(x_{n}-\cdot)\}=n
$$
for every choice of $0\leq x_{1}<\cdots<x_{n}<2\pi$ and all $n\in
{\Bbb N}$.
\end{definition}

 Denote by  $L_q$:= $L_q[0, 2\pi]$  the classical
Lebesgue integral space of 2$\pi$--periodic real--valued functions
on ${\mathbb R}$ with the usual norm $\|\cdot\|_q$, $1\leq
q\leq\infty$. We now introduce the classes of functions to be
studied here. Let $\varphi$ be a differentiable, odd and strictly
increasing function for which ${\varphi}'$ is continuous on
$[-1,1]$. It is clear that the functions
\begin{equation}
\varphi_{0}(z)=\tan(\pi z/4) \label{dt4-7}
\end{equation}
and
\begin{equation}
\varphi_{1}(z)=z \label{dt4-8}
\end{equation}
satisfy the above conditions. Let $f\ast g$ denote the convolution
of the functions $f$ and $g$, i.e.,
$$
(f\ast g)(x):=\frac{1}{2\pi}\int_0^{2\pi}f(x-t)g(t)\,dt,
$$
and let $G$ be a kernel satisfying  Property $B$ or $NCVD$. Then $G$
has some variation-diminishing properties, we are interested in the
class of $2\pi$--periodic functions representable as
 as follows:
\begin{equation}
\widetilde{K}^{G,\varphi}_{\infty,\beta}=\{f:\,
f=a+G\ast\varphi(K_{\beta}\ast u),\, a\in\Theta,\,
\varphi(K_{\beta}\ast u)\perp \Theta,\,\|u\|_{\infty} \leq 1\},
\label{dt4-c}
\end{equation}
where $\Theta={\mathbb R}$ if $G$ satisfies Property $B$ or
 $\Theta=\emptyset$ if $G$ is a $NCVD$ kernel, and
\begin{equation}
K_{\beta}(z):=1+2\sum_{k=1}^{\infty}\frac{\cos(kz)}{\cosh(k\beta)},
\qquad \beta >0. \label{K}
\end{equation}
For conciseness, we will regard $\varphi(K_{\beta}\ast u)$ as $u$
if $\beta=0$. It can be seen that
$\widetilde{K}^{G,\varphi_{0}}_{\infty,\beta}$ and
$\widetilde{K}^{G,\varphi_{1}}_{\infty,\beta}$ are convex. Set
\begin{equation}
D_{r}(t):=2\sum_{k=1}^{\infty}\frac{\cos(kt-(\pi
r/2))}{k^{r}},\qquad r=1,2,\ldots.\label{D}
\end{equation}

It is known that for each $r\geq 2$, $D_{r}$ satisfies Property
$B$ ( $D_{1}$ also satisfies all the conditions of Property $B$
except that it has a jump discontinuity at $x=0$ )( see \cite{P},
p.133). Then from \cite{Oc97b} we have
\begin{equation}
\widetilde{K}^{G,\varphi}_{\infty,\beta}:=\left\{
\begin{array}{ll}
\widetilde{W}^{r}_{\infty},\quad &G=D_{r},\,\beta=0,\\[2mm]
\widetilde{h}^{r}_{\infty,\beta},\quad &G=D_{r},\,\varphi=\varphi_{1},\,\beta>0,\\[2mm]
\widetilde{H}^{r}_{\infty,\beta},\quad
&G=D_{r},\,\varphi=\varphi_{0},\,\beta>0,
\end{array}
\right.\label{dt4-9}
\end{equation}
where $\widetilde{W}^{r}_{\infty}$ is the classical Sobolev class
of real, 2$\pi$-periodic functions $f$ whose $(r-1)$th derivative
is absolutely continuous and whose $r$th derivative satisfies the
condition $\|f^{(r)}\|_\infty\leq 1$;
$\widetilde{h}^{r}_{\infty,\beta}$ is the Achieser class
\cite[p. 214, and p. 219]{Ac56}
 which are all real--valued on ${\mathbb R}$ with $2\pi$--periodic,
and can be continued analytically in the
strip $S_{\beta}:=\{z\in {\mathbb C}:|{\rm
Im}z|<\beta\}$,  satisfying  the restriction condition $|{\rm Re}
f^{(r)}(z)|\leq 1$ in this strip;
and $\widetilde{H}^{r}_{\infty,\beta}$ is the class of
$2\pi$--periodic, real--valued functions on ${\mathbb R}$ which
are analytic in the strip $S_{\beta}$ and satisfy the condition
$|f^{(r)}(z)|\leq 1$, $z\in S_{\beta}$. Thus the class
$\widetilde{K}^{G,\varphi}_{\infty,\beta}$ includes the three famous
classes as its special cases. \\

The exact values of $n$--widths were investigated for
the classical Sobolev class $\widetilde{W}^{r}_{p}$ (see
\cite{Ko}) in $L_q$ space, by the efforts of many mathematicians
some similar results were also built for the classes
$\widetilde{h}^{r}_{\infty,\beta}$ and
$\widetilde{H}^{r}_{\infty,\beta}$. The Kolmogorov $n$--widths of
the class $\widetilde{h}^{r}_{\infty,\beta}$ in the space
$L_\infty$ were obtained by Tikhomirov \cite{T}, Forst \cite{Fo}
for $r=0$ and by Sun \cite{S}, and Osipenko \cite{Oc97a} for all
$r\in{\mathbb N}$. After this, the exact values of the even
$n$--widths of the class of $\widetilde{h}^{r}_{\infty,\beta}$ in
$L_q$, $1\leq q<\infty$ were calculated by Pinkus \cite{P} for
$r=0$ and by Osipenko \cite{Oc95b} for all $r\in{\mathbb N}$. The
exact estimates of the even $n$--widths of the class
$\widetilde{H}^{r}_{\infty,\beta}$ in the space $L_q$ were
determined by Osipenko \cite{Oc94}  for $r=0$, $1\leq q\leq\infty$
and by Osipenko \cite{Oc97a}--\cite{Oc97b}, Osipenko and
Wilderotter \cite{Ow} for all $r\in{\mathbb N}$, $q=\infty$. The
exact values of Gel$'$fand $n$--width, Kolmogorov $2n$--width,
linear $2n$--width and information $2n$--width  of
$\widetilde{H}^r_{\infty,\beta}$ in the space $L_q$ were obtained
by us \cite{FGSLXH}  for all
$r\in{\mathbb N}$, $1\leq q<\infty$ .

 Now we outline the rest of this paper. In Section \ref{sec2},
we establish a comparison theorem of Kolmogorov type on
$\widetilde{K}^{G,\varphi}_{\infty,\beta}$ by using the property
of the kernel $G$. In Section \ref{sec3}, using the results of Section \ref{sec2},
we get a Taikov type inequality which will be used as the upper
estimates of the Gel'fand $n$-width of the class
$\widetilde{K}^{G,\varphi}_{\infty,\beta}$ in $L_{q}$ for $G$
satisfying Property $B$ and $1\leq q<\infty$. In Section \ref{sec4},
following  the method of \cite{FGSLXH} and \cite{Oc97b}, we solve
two minimum norm questions on analogue of the classical polynomial
perfect splines, and then we use them to prove the lower estimates
of $n$-widths, which together with some results of Section \ref{sec3}
determine the precise values of $n$--widths of
$\widetilde{K}^{G,\varphi}_{\infty,\beta}$ in the space
$L_{\infty}$ for $\varphi$ being such that
$\widetilde{K}_{\infty,\beta}^{G,\varphi}$ is convex, the exact
values of the Gel'fand $n$--width, and the lower estimate of
the Kolmogorov $2n$--width of
$\widetilde{K}^{G,\varphi}_{\infty,\beta}$ in $L_{q}$ for $G$
satisfying Property $B$ and $1\leq q<\infty$.

\section{ Comparison Theorem of Kolmogorov Type}\label{sec2}

\hspace{1ex}  The Kolmogorov comparison theorem (see \cite{Ko}),
which concerns the comparison of derivatives of  differentiable
functions defined on the real line, plays an important role in
establishing some sharp estimates of extremal problems
in approximation theory. In
this section, we will prove a comparison theorem of Kolmogorov
type on the class $\widetilde{K}^{G,\varphi}_{\infty,\beta}$.
 Before advancing our discuss further, we introduce
the following notions. Let
\begin{equation}
\Lambda_{2n}=\{\xi:\xi=(\xi_{1},\ldots,\xi_{2n}),\,
0\leq\xi_{1}<\cdots<\xi_{2n}<2\pi\},\quad n\in {\mathbb N}.
\label{g1}
\end{equation}
The closure of $\Lambda_{2n}$, $\overline{\Lambda}_{2n}$ is given
by
\begin{equation}
\overline{\Lambda}_{2n}:=\{\xi:\xi=(\xi_{1},\ldots,\xi_{2m}),\,0\leq
\xi_{1}<\cdots<\xi_{2m}<2\pi,\ m\leq n\}. \label{g2}
\end{equation}
For each $\xi\in{\Lambda}_{2n}$, we define
\begin{equation} h_{\xi}(t):=(-1)^{j},\quad
t\in[\xi_{j-1},\xi_{j}),\quad j=1,\ldots,2n+1, \label{g3}
\end{equation}
where $\xi_{0}=0$ and $\xi_{2n+1}=2\pi$. There is one particularly
important function of this form, say $h_n$, which corresponds to
the choice $\xi_{j}=(j-1)\pi/n$, $j=1,\ldots,2n$.  \\

 We introduce the standard function
$\Phi_{n,\beta}^{G,\varphi}$ of the comparison theorem on the
class $\widetilde{K}^{G,\varphi}_{\infty,\beta}$ and discuss some
of its properties. Let
\begin{equation}
\Phi_{n,\beta}^{G,\varphi}(x):=(G\ast\varphi(K_{\beta}\ast
h_{n}))(x).\label{dt4-10}
\end{equation}
Then it is easily seen that
$\Phi_{n,\beta}^{G,\varphi}(x+\pi/n)=-\Phi_{n,\beta}^{G,\varphi}(x)$
and $\Phi_{n,\beta}^{G,\varphi}(x)$ is continuous with respect to
$x$. So $\Phi_{n,\beta}^{G,\varphi}$ is $2\pi/n$--periodic and
there exist at least $2n$ points $0\leq t_{1}<\cdots<t_{2n}<2\pi$
such that
\begin{equation}
\Phi_{n,\beta}^{G,\varphi}(t_{j}):=\epsilon
(-1)^{j}\|\Phi_{n,\beta}^{G,\varphi}\|_{\infty},\quad
j=1,\ldots,2n, \label{dt4-11}
\end{equation}
where $\epsilon=1$ or $ -1$.\\

 Now we are in position to prove the comparison
theorem on
the class $\widetilde{K}^{G,\varphi}_{\infty,\beta}$.

\begin{theorem}[Comparison theorem of Kolmogorov type]\label{Theorem dt4-1}
Let $f\in\widetilde{K}^{G,\varphi}_{\infty,\beta}$ be such that
$\|f\|_\infty\leq \|\Phi_{n,\beta}^{G,\varphi}\|_\infty$ for some
positive integer $n$, and
$f(\alpha)=\Phi_{n,\beta}^{G,\varphi}(\gamma)$,
$f'(\alpha){\Phi_{n,\beta}^{G,\varphi}}'(\gamma)\geq 0$ for some
$\alpha,\,\gamma\in {\mathbb R}$. Then
\begin{equation}
|f'(\alpha)|\leq|{\Phi_{n,\beta}^{G,\varphi}}'(\gamma)|.
\label{dt4-12}
\end{equation}
\end{theorem}
\begin{proof} Without loss of generality
we may assume $\alpha=\gamma$.
 Assume that for some $f=a+G\ast\varphi(K_{\beta}\ast u)
 \in\widetilde{K}^{G,\varphi}_{\infty,\beta}$, $a\in\Theta$,
$\varphi(K_{\beta}\ast u)\perp \Theta$, $\|u\|_{\infty} \leq 1$,
positive integer $n$ and $\alpha\in{\mathbb R}$, we have
$\|f\|_\infty\leq\|\Phi_{n,\beta}^{G,\varphi}\|_\infty$,
$f(\alpha)=\Phi_{n,\beta}^{G,\varphi}(\alpha)$,
$f'(\alpha){\Phi_{n,\beta}^{G,\varphi}}'(\alpha)\geq 0$ and
$|f'(\alpha)|>|{\Phi_{n,\beta}^{G,\varphi}}'(\alpha)|$. It follows
from the definitions  that $f'$ and
${\Phi_{n,\beta}^{G,\varphi}}'$ are continuous. If $f(\alpha)\neq
0$, the continuity of them ensures the existence of $0<\rho<1$ and
$\alpha_{0}\in {\mathbb R}$ such that
\begin{equation}
\rho f(\alpha_{0})=\Phi_{n,\beta}^{G,\varphi}(\alpha_{0}),\quad
\rho|f'(\alpha_{0})|>|{\Phi_{n,\beta}^{G,\varphi}}'(\alpha_{0})|.\label{dt4-13}
\end{equation}
If $f(\alpha)=\Phi_{n,\beta}^{G,\varphi}(\alpha)=0$, we can choose
some $\rho$ satisfying $1> \rho >
\frac{|{\Phi_{n,\beta}^{G,\varphi}}'(\alpha)|}{|f'(\alpha)|}$\,
and let $\alpha_{0}=\alpha$. Then (\ref{dt4-13}) also holds. Put
$\overline{f}:=\rho f$, then
$\overline{f}\in\widetilde{K}^{G,\varphi}_{\infty,\beta}$ and
$\|\overline{f}\|_\infty< \|\Phi_{n,\beta}^{G,\varphi}\|_\infty$.
It is enough to consider one possible case
$\overline{f}(\alpha_{0})=\Phi_{n,\beta}^{G,\varphi}(\alpha_{0})\geq
0$,
${\overline{f}}'(\alpha_{0})>{\Phi_{n,\beta}^{G,\varphi}}'(\alpha_{0})>0$
(the other cases can be treated by the same way). From the
property of $\Phi_{n,\beta}^{G,\varphi}$, it has at least $2n$
monotone intervals on every interval whose length is $2\pi$.
Choose a interval whose length is $2\pi$ and whose endpoints are
extremal points of $\Phi_{n,\beta}^{G,\varphi}$, in which
$\alpha_{0}$ is contained, written as $\Delta_{\alpha_{0}}$. By
geometrical consideration we see that in the monotone interval,
which is in $\Delta_{\alpha_{0}}$ and contains $\alpha_{0}$, the
graphs of $\overline{f}$ and $\Phi_{n,\beta}^{G,\varphi}$
intersect at least three times, while on each of the other
monotone intervals in $\Delta_{\alpha_{0}}$ these graphs intersect
at least once. Hence for $F(x):=\Phi_{n,\beta}^{G,\varphi}(x)-\rho
f(x)$, we have
$\widetilde{Z}_{c}(F)\geq 2n+2$.

 On the other hand, since $\varphi$ is a
differentiable odd and strictly increasing function, and
$K_{\beta}$ is $NCVD$, if $G$ satisfies Property $B$, by
\cite[Chapt. IV, Prop. 6.4]{P} we can get
$$
\widetilde{Z}_c(F)\leq S_{c}(\varphi(K_{\beta}\ast
h_{n})-\rho\varphi(K_{\beta}\ast u))\leq S_{c}(h_{n}(\cdot)-
u^{\star}(\cdot))\leq 2n,
$$
where $u^{\star}$ is defined by the equality
$\rho\varphi(K_{\beta}\ast u)=\varphi(K_{\beta}\ast u^{\star})$,
satisfying $\varphi(K_{\beta}\ast u^{\star})\perp \Theta$ and
$\|u^{\star}\|_{\infty}<1$. If $G$ is a $NCVD$ kernel, by
\cite[Chapt. III, Section 3]{P}, we shall assume that $G$ is
extended $CVD$ by means of ``smoothing". The property of $G$ being
extended $CVD$ implies that
$$
\widetilde{Z}_c(F)\leq S_{c}(\varphi(K_{\beta}\ast
h_{n})-\rho\varphi(K_{\beta}\ast u))\leq S_{c}(h_{n}(\cdot)-
u^{\star}(\cdot))\leq 2n,
$$
where $u^{\star}$ is the same as that above. Thus we get a
contradiction. Theorem \ref{Theorem dt4-1} is proved.
\end{proof}

\begin{remark}\label{remark dt4-1} Among others, Sun \cite{S} established
a comparison theorem of Kolmogorov type for the class
$\widetilde{h}_{\infty, \beta}$. Moreover, Fisher \cite{Fi}
demonstrated an inequality of Landau--Kolmogorov type for the
class of non--periodic analytic functions on the open unit disk
whose $r$th derivatives are bounded by $1$ and which  are
real--valued on the interval $(-1, 1)$. Osipenko \cite{Oc94b}
established Kolmogorov inequalities for the class of functions
analytic in the strip $S_{\beta}$, real on the real axis
$\mathbb{R}$ and satisfying the condition $|f^{(r)}(z)|\leq 1$,
$z\in S_{\beta}$, which includes the class $\widetilde{H}_{\infty, \beta}^{r}$
as a subset with 2$\pi$-periodic.
Using a method different from  that of Osipenko
\cite{Oc94b} we established a comparison theorem of Kolmogorov
type on  the class $\widetilde{H}_{\infty, \beta}^{r}$ (see
\cite{FGSLXH}).
\end{remark}

 Our next aim is to establish a Landau--Kolmogorov
type inequality and some other corollaries of the comparison
theorem.

\begin{corollary}[Landau--Kolmogorov type inequality]\label{Corollary dt4-1} Let
$f\in\widetilde{K}^{G,\varphi}_{\infty,\beta}$ be such that
$\|f\|_\infty\leq \|\Phi_{n,\beta}^{G,\varphi}\|_\infty$ for some
positive integer $n$ . Then
\begin{equation}
\|f'\|_{\infty}\leq \|{\Phi_{n,\beta}^{G,\varphi}}'\|_\infty .
\label{dt4-f1}
\end{equation}
\end{corollary}
\begin{proof} Assume that there exists a
$\xi \in {\mathbb R}$ for which $|f'(\xi)|>
\|{\Phi_{n,\beta}^{G,\varphi}}'\|_\infty$. Since
$\Phi_{n,\beta}^{G,\varphi}(t)$ is a continuous function of $t$
and
$$
\|f\|_\infty\leq \|\Phi_{n,\beta}^{G,\varphi}\|_\infty,\qquad
-\|\Phi_{n,\beta}^{G,\varphi}\|_\infty\leq
\Phi_{n,\beta}^{G,\varphi}(t)\leq
\|\Phi_{n,\beta}^{G,\varphi}\|_\infty,\quad t\in {\mathbb R},
$$
there exists an $\eta \in {\mathbb R}$ such that  $f(\xi)=
\Phi_{n,\beta}^{G,\varphi}(\eta)$ and $|f'(\xi)|>
|{\Phi_{n,\beta}^{G,\varphi}}'(\eta)|$, which contradicts Theorem
\ref{Theorem dt4-1}. Corollary \ref{Corollary dt4-1} is
proved.
\end{proof}

 Let $C({\mathbb R})$ be the set of all continuous
functions on the real line ${\mathbb R}$. A function $\psi\in
C({\mathbb R})$ is said to be regular (see \cite{Ko}, p. 107) if
it has a period $2l$ and if on some interval $(a,a+2l)$ ( $a$ is a
point of an absolute extremum of $\psi$ ) there is a point $c$
such that $\psi$ is strictly monotone on $(a,c)$ and $(c,a+2l)$.
In order to emphasize the length of the period $2l$,
sometimes we shall speak about $\psi$ as being $2l$--regular.

 We say that a function $f\in C({\mathbb R})$
possesses a $\mu$--property with respect to a regular function
$\psi$ if for every $\alpha\in {\mathbb R}$ and on every interval
of monotonicity of $\psi$ the difference $\psi(t)-f(t+\alpha)$
either does not change sign or changes sign exactly once--from $+$
to $-$ if $\psi$ decreases or from $-$ to $+$ if $\psi$ increases.
It is clear that $f$ will possess the $\mu$--property with respect
to $\psi (t+\beta)$, $\beta\in {\mathbb R}$ if it
possesses the $\mu$--property with respect to $\psi$.

\begin{corollary}\label{Corollary dt4-2} Let
$G$ satisfy Property $B$. Let
$f\in\widetilde{K}^{G,\varphi}_{\infty,\beta}$ be such that
$\|f\|_\infty\leq \|\Phi_{n,\beta}^{G,\varphi}\|_\infty$ for some
positive integer $n$. Then the function $f$ possesses a
$\mu$--property with respect to $\Phi_{n,\beta}^{G,\varphi}$.
\end{corollary}
\begin{proof}
 Let $G$ satisfy Property
$B$. First we will prove that $\Phi_{n,\beta}^{G,\varphi}$ is
$2\pi/n$--regular. Assume that $\Phi_{n,\beta}^{G,\varphi}$ is not
$2\pi/n$--regular. Since
$\Phi_{n,\beta}^{G,\varphi}(x+\pi/n)=-\Phi_{n,\beta}^{G,\varphi}(x)$
and $\Phi_{n,\beta}^{G,\varphi}(x)$ is continuous with respect to
$x$, there exists a constant $c\in {\mathbb R}$ such that
$\Phi_{n,\beta}^{G,\varphi}(\cdot)-c$ has at least $2n+2$
zeros, i.e., $\widetilde{Z}_c(\Phi_{n,\beta}^{G,\varphi}(\cdot)-c)\geq 2n+2$.

 On the other hand, since $\varphi$ is a
differentiable, odd and strictly increasing function, and
$K_{\beta}$ is $NCVD$, the condition that $G$ satisfies Property
$B$ implies that
$$
\widetilde{Z}_c(\Phi_{n,\beta}^{G,\varphi}(\cdot)-c)\leq
S_{c}(\varphi(K_{\beta}\ast h_{n}))\leq S_{c}(h_{n})\leq 2n
$$
by \cite[Chapt. IV, Prop. 6.4]{P}. This contradiction shows that
$\Phi_{n,\beta}^{G,\varphi}$ is $2\pi/n$--regular. Then there
exist exactly $2n$ points $0\leq
t_{n,\beta,1}^{G,\varphi}<\cdots<t_{n,\beta,2n}^{G,\varphi}<2\pi$
such that
\begin{equation}
\Phi_{n,\beta}^{G,\varphi}(t_{n,\beta,j}^{G,\varphi}):=\epsilon
(-1)^{j}\|\Phi_{n,\beta}^{G,\varphi}\|_{\infty},\quad
t_{n,\beta,j+1}^{G,\varphi}-t_{n,\beta,j}^{G,\varphi}=\frac{\pi}{n},\quad
j=1,\ldots,2n, \label{ff3}
\end{equation}
where $\epsilon=1$ or $ -1$, and
$t_{n,\beta,2n+1}^{G,\varphi}=t_{n,\beta,1}^{G,\varphi}+2\pi$. So
put
\begin{equation}
\Delta_{n,\beta}^{G,\varphi,(j)}=[t_{n,\beta,j}^{G,\varphi},t_{n,\beta,j+1}^{G,\varphi}),\quad
j=1,\ldots,2n. \label{dt4-11r}
\end{equation}
Then $\Phi_{n,\beta}^{G,\varphi}$ is monotonic on each interval
$\Delta_{n,\beta}^{G,\varphi,(j)}$, $j=1,\ldots,2n$.

Suppose that $f$ does not possess a $\mu$--property
with respect to $\Phi_{n,\beta}^{G,\varphi}$. Then for some
$\alpha\in {\mathbb R}$, the difference
$\Phi_{n,\beta}^{G,\varphi}(t)-f(t+\alpha)$ changes sign on a
monotone interval $\Delta_{n,\beta}^{G,\varphi,(k)}$
($k\in\{1,\ldots,2n\}$) of the standard function
$\Phi_{n,\beta}^{G,\varphi}$ from $+$ to $-$ if
$\Phi_{n,\beta}^{G,\varphi}$ increases and from $-$ to $+$ if
$\Phi_{n,\beta}^{G,\varphi}$ decreases. By continuity arguments
this fact will also hold for the difference
$$
F(t):=\Phi_{n,\beta}^{G,\varphi}(t)-(1-\epsilon)f(t+\alpha)
$$
for $0<\epsilon<1$ sufficiently small. Since
$(1-\epsilon)\|f(\cdot+\alpha)\|_{\infty}<\|\Phi_{n,\beta}^{G,\varphi}\|_\infty$,
$F$ changes sign at least three times on the interval
$\Delta_{n,\beta}^{G,\varphi,(k)}$ and at least once on each of
the other intervals $\Delta_{n,\beta}^{G,\varphi,(j)}$
$(j\in\{1,\ldots,2n\}\backslash\{k\})$. Hence $S_{c}(F)\geq
2n+2$.

 On the other hand, we can prove that $S_{c}(F)\leq
2n$ by the same method as that in the proof of Theorem
\ref{Theorem dt4-1}. This contradiction proves the
corollary.
\end{proof}
 A function $F$ is said to be a periodic integral of a
real, $2\pi$-periodic, continuous function $f$ on ${\mathbb R}$ if
$F'(x)=f(x)$ and $F(x+2\pi)=F(x)$ for all $x\in {\mathbb R}$.
 Let $G$ satisfy Property $B$. Put
\begin{equation}
\widetilde{G}(x)=\int_{0}^{x}(G(y)-a_{0})\,dy, \label{dt4-G1}
\end{equation}
where $a_{0}=1/2\pi\int_{0}^{2\pi}G(y)dy$. Then $\widetilde{G}$
satisfies Property $B$ (see \cite[Chapt. IV, Prop. 6.6]{P}) and
$\Phi_{n,\beta}^{\widetilde{G},\varphi}$ is a periodic integral of
$\Phi_{n,\beta}^{G,\varphi}$. From Corollary \ref{Corollary
dt4-1}, we have immediately the following
theorem which will be used in the next section.

\begin{theorem}\label{Theorem dt4-2} Let $G$ satisfy Property $B$ and
$\widetilde{G}$ be defined by (\ref{dt4-G1}). Suppose that
 $f\in \widetilde{K}^{G,\varphi}_{\infty,\beta}$, and $F$ is a periodic
integral of $f$ such that $\|F\|_\infty\leq
\|\Phi_{n,\beta}^{\widetilde{G},\varphi}\|_\infty$ for some
positive integer $n$. Then
\begin{equation}
 \|f\|_{\infty}\leq\|\Phi_{n,\beta}^{G,\varphi}\|_{\infty}.
\label{dt4-f8}
\end{equation}
\end{theorem}

\section{Inequalities of Taikov Type}\label{sec3}

\hspace{1ex} In this section, we establish an inequality of Taikov type,
which leads to the upper estimates of the Gel'fand $n$-widths of
$\widetilde{K}^{G,\varphi}_{\infty,\beta}$ in $L_{q}$ for $G$
satisfying Property $B$ and $1\leq q<\infty$. To do this, we need
some auxiliary lemmas.

\begin{lemma}\label{Lemma dt4-1}  Let $n\in {\mathbb N}$, $G$ satisfy
Property $B$ and $\widetilde{G}$ be defined by (\ref{dt4-G1}).
Then
\begin{equation}
\int_{0}^{2\pi}|\Phi_{n,\beta}^{G,\varphi}(t)|\,dt
=\bigvee_{0}^{2\pi}\Phi_{n,\beta}^{\widetilde{G},\varphi}
=4n\|\Phi_{n,\beta}^{\widetilde{G},\varphi}\|_{\infty}.
\label{dt4-r}
\end{equation}
\end{lemma}

\begin{proof}\ By (\ref{dt4-11}) and
(\ref{dt4-11r}), $\Phi_{n,\beta}^{G,\varphi}$ is
$2\pi/n$-periodic,  strictly monotonic on
$\Delta_{n,\beta}^{G,\varphi,(j)}$ ($j=1,\ldots,2n$) and
$\Phi_{n,\beta}^{G,\varphi}(t+\frac{\pi}{n})=-\Phi_{n,\beta}^{G,\varphi}(t)$
for all $t\in[0,2\pi)$. Hence,
$$
\begin{array}{lll}
&\dint_{0}^{2\pi}|\Phi_{n,\beta}^{G,\varphi}(t)|\,dt
=\dint_{t_{n,\beta,1}^{\widetilde{G},\varphi}}^{t_{n,\beta,1}^{\widetilde{G},\varphi}
+2\pi}|\Phi_{n,\beta}^{G,\varphi}(t)|\,dt
=2n\dint_{t_{n,\beta,1}^{\widetilde{G},\varphi}}
^{t_{n,\beta,1}^{\widetilde{G},\varphi}+
\frac{\pi}{n}}|\Phi_{n,\beta}^{G,\varphi}(t)|\,dt=2n\Bigg|\dint_{t_{n,\beta,1}^{\widetilde{G},\varphi}}
^{t_{n,\beta,1}^{\widetilde{G},\varphi}+
\frac{\pi}{n}}\Phi_{n,\beta}^{G,\varphi}(t)\,dt\Bigg|\\
=&2n\left|\Phi_{n,\beta}^{\widetilde{G},\varphi}(t_{n,\beta,1}^{\widetilde{G},\varphi}+
\frac{\pi}{n})-\Phi_{n,\beta}^{\widetilde{G},\varphi}
(t_{n,\beta,1}^{\widetilde{G},\varphi})\right|=4n\|\Phi_{n,\beta}^{\widetilde{G},\varphi}\|_{\infty}
=2n\dbigvee_{t_{n,\beta,1}^{\widetilde{G},\varphi}}^{t_{n,\beta,1}^{\widetilde{G},\varphi}+
\frac{\pi}{n}}\Phi_{n,\beta}^{\widetilde{G},\varphi}
=\dbigvee_{0}^{2\pi}\Phi_{n,\beta}^{\widetilde{G},\varphi}.
\end{array}
$$
\end{proof}

 Now let $f$ be a $2\pi$-periodic integrable function,
and denote by $r(f,t)$ the non-increasing rearrangement of $|f|$
(see \cite {Ko}, p. 110). With this notation we have

\begin{lemma}\label{Lemma dt4-2} (see \cite{Ko}, p.112)\ \ Let $f$, $g\in L_q$
$(1\leq q<\infty)$ and
$$
\int_{0}^{x}r(f,t)\,dt \leq\int_{0}^{x}r(g,t)\,dt,\qquad 0\leq
x\leq 2\pi.
$$
Then
$$
\|f\|_q\leq\|g\|_q.
$$
\end{lemma}

\begin{lemma}\label{Lemma dt4-3} (see \cite{Ko}, p. 114) Let $f\in C({\mathbb R)}$
possess the $\mu$-property with respect to the $2\pi/n$--regular
($n=1,2,\ldots$) function $\psi$ and
$$
\int_0^{2\pi/n}\psi(t)\,dt=0.
$$
If
$$
\ \ \ \ \ \ \ \dmin_{u}\psi(u)\leq f(t)\leq\dmax_u\psi(u), \qquad
\forall\,t\in{\mathbb R},
$$
and
$$
\dmax_{a,b}\left|\dint_a^bf(t)\,dt\right|\leq\frac 1
2\dint_0^{2\pi/n}|\psi(t)|\,dt,
$$
then
$$
\int_{0}^{x}r(f,t)\,dt \leq\int_{0}^{x}r(\psi,t)\,dt,\qquad 0\leq
x\leq 2\pi.
$$
\end{lemma}

\begin{theorem}\label{Theorem dt4-3}  Let $G$ satisfy
Property $B$, $\widetilde{G}$ be defined by (\ref{dt4-G1}),
$f\in\widetilde{K}^{G,\varphi}_{\infty,\beta}$ and $F$ be a
periodic integral of $f$ such that $\|F\|_\infty\leq
\|\Phi_{n,\beta}^{\widetilde{G},\varphi}\|_\infty$ for some
positive integer $n$. Then
\begin{equation}
\int_{0}^{x}r(f,t)\,dt
\leq\int_{0}^{x}r(\Phi_{n,\beta}^{G,\varphi},t)\,dt,\quad 0\leq
x\leq 2\pi, \label{dt4-g1}
\end{equation}
\begin{equation}
\|f\|_q \leq\|\Phi_{n,\beta}^{G,\varphi}\|_q, \quad 1\leq
q<\infty. \label{dt4-g2}
\end{equation}
\end{theorem}
\begin{proof} By virtue of Lemma
\ref{Lemma dt4-2}, the inequality (\ref{dt4-g2}) follows from
(\ref{dt4-g1}). So we only need to prove the inequality
(\ref{dt4-g1}). Since
$f\in\widetilde{K}^{G,\varphi}_{\infty,\beta}$ and $F$ is a
periodic integral of $f$ such that
$\|F\|_{\infty}\leq\|\Phi_{n,\beta}^{\widetilde{G},\varphi}\|_\infty$
for some positive integer $n$,  we conclude from  Theorem
\ref{Theorem dt4-2} that $\|f\|_\infty\leq
\|\Phi_{n,\beta}^{G,\varphi}\|_\infty$, and it follows from
Corollary \ref{Corollary dt4-2} that $f$ possesses a
$\mu$-property with respect to the $2\pi/n$--regular function
$\Phi_{n,\beta}^{G,\varphi}$. By noticing that
$\int_{0}^{2\pi/n}\Phi_{n,\beta}^{G,\varphi}(t) \,dt=0$ and the
following result derived from the proof of Lemma \ref{Lemma
dt4-1}:
$$
\begin{array}{ll}
&\ \ \dmax_{a,b}|\dint_{a}^{b} f(t) dt|\leq 2\|F\|_{\infty}
\leq 2\|\Phi_{n,\beta}^{\widetilde{G},\varphi}\|_\infty\\
&=\dint_{0}^{\pi/n}
|\Phi_{n,\beta}^{G,\varphi}(t)|\,dt=\frac{1}{2}\dint_{0}^{2\pi/n}
|\Phi_{n,\beta}^{G,\varphi}(t)|\,dt,
\end{array}
$$
we see that Lemma \ref{Lemma dt4-3} is applicable for the
functions $f$ and $\Phi_{n,\beta}^{G,\varphi}$. Hence the
inequality (\ref{dt4-g1}) is true. Theorem \ref{Theorem dt4-3} is
proved.
\end{proof}

 Let
$$
\mathcal{T}_{n}={\rm span}\{1, \cos t,\sin
t,\ldots,\cos((n-1)t),\sin((n-1)t)\}
$$
be the trigonometric polynomial subspace of order $n-1$. Set
$$
a_{j}(f):=\frac{1}{\pi}\int_{0}^{2\pi}f(t)\cos(jt)\,dt,\quad
j=0,1,\ldots,
$$
$$
b_{j}(f):=\frac{1}{\pi}\int_{0}^{2\pi}f(t)\sin(jt)\,dt,\quad
j=1,2,\ldots,
$$
$$
I_{2n-1}(f):=(a_{0}(f),a_{1}(f),b_{1}(f),\cdots,a_{n-1}(f),b_{n-1}(f)).
$$
The following lemma will be used in establishing Taikov type inequality in the next
section as well as the upper estimates of $n$-widths in section \ref{sec4}.
\begin{lemma}\label{Lemma dt4-4}\ Let $n\in {\mathbb N}$ and  $\varphi$
be such that $\widetilde{K}^{G,\varphi}_{\infty,\beta}$ is convex.
Then there exists a linear operator ${\mathcal L}$: ${\mathbb
R}^{2n-1}\rightarrow \mathcal{T}_{n}$ such that
\begin{equation}
 \dsup_{f\in\widetilde{K}^{G,\varphi}_{\infty,\beta}}\|f-{\mathcal L}(I_{2n-1}(f))\|_{\infty}
 =\|G\ast \varphi(K_{\beta}\ast
h_{n})\|_{\infty}=\|\Phi_{n,\beta}^{G,\varphi}\|_{\infty}.
 \label{dt4-18}
\end{equation}
\end{lemma}
\begin{proof} First, we claim that
\begin{equation}
\dsup_{f\in\widetilde{K}^{G,\varphi}_{\infty,\beta},
I_{2n-1}(f)=0}\|f\|_{\infty}
 =\|G\ast \varphi(K_{\beta}\ast
h_{n})\|_{\infty}=\|\Phi_{n,\beta}^{G,\varphi}\|_{\infty}.
\label{dt4-19}
\end{equation}
Suppose that there exists a function $f_{0}=a_{0}+G\ast
\varphi(K_{\beta}\ast
u_{0})\in\widetilde{K}^{G,\varphi}_{\infty,\beta}$, where
$a_{0}\in\Theta$, $\varphi(K_{\beta}\ast u_{0})\perp \Theta$ and
$\|u_{0}\|_{\infty} \leq 1$, such that $I_{2n-1}(f_{0})=0$ and
$\|f_{0}\|_{\infty}>\|\Phi_{n,\beta}^{G,\varphi}\|_{\infty}$. Set
$$
\rho=\frac{\|\Phi_{n,\beta}^{G,\varphi}\|_{\infty}}{\|f_{0}\|_{\infty}}.
$$
Then $0<\rho<1$. We choose $\alpha\in [0,2\pi)$ such that
$F(\cdot)\triangleq \Phi_{n,\beta}^{G,\varphi}(\cdot+\alpha)-\rho
f_{0}(\cdot)$ has a multiple zero. Since
$\Phi_{n,\beta}^{G,\varphi}$ has the period $2\pi/n$,
$I_{2n-1}(\Phi_{n,\beta}^{G,\varphi})=0$. Thus, $I_{2n-1}(F)=0$.
Since the trigonometric system is a Tchebycheff system, it follows
from \cite[p.41]{P} that
$$
\widetilde{Z}_c(F(\cdot))>S_{c}(F(\cdot))\geq 2n.
$$
\hspace{1ex} On the other hand, since $\varphi$ is a
differentiable odd and strictly increasing function, and
$K_{\beta}$ is $NCVD$, if $G$ satisfies Property $B$, by
\cite[Chapt. IV, Prop. 6.4]{P} we can get
$$
\begin{array}{ll}
\widetilde{Z}_c(F(\cdot))&\leq S_{c}(\varphi(K_{\beta}\ast
h_{n})(\cdot+\alpha)-\rho\varphi(K_{\beta}\ast u_{0})(\cdot))\\
&\leq S_{c}(h_{n}(\cdot+\alpha)- u_{0}^{\star}(\cdot))\leq 2n,
\end{array}
$$
where $u_{0}^{\star}$ is defined by the equality
$\rho\varphi(K_{\beta}\ast u_{0})=\varphi(K_{\beta}\ast
u_{0}^{\star})$, satisfying $\varphi(K_{\beta}\ast
u_{0}^{\star})\perp \Theta$ and $\|u_{0}^{\star}\|_{\infty}<1$. If
$G$ is a $NCVD$ kernel, by \cite[Chapt. III, Section 3]{P}, we
shall assume that $G$ is extended $CVD$ by means of ``smoothing".
The property of $G$ being extended $CVD$ implies that
$$
\begin{array}{ll}
\widetilde{Z}_c(F(\cdot))&\leq S_{c}(\varphi(K_{\beta}\ast
h_{n})(\cdot+\alpha)-\rho\varphi(K_{\beta}\ast u_{0})(\cdot))\\
&\leq S_{c}(h_{n}(\cdot+\alpha)- u_{0}^{\star}(\cdot))\leq 2n,
\end{array}
$$
where $u_{0}^{\star}$ is the same as that above. Thus we get a
contradiction. So (\ref{dt4-19}) holds.

 Now we consider the problem of the optimal recovery
of the value $f(0)$ in the class
$\widetilde{K}^{G,\varphi}_{\infty,\beta}$ from the information
$I_{2n-1}(f)$. Since $\varphi$ is such that
$\widetilde{K}^{G,\varphi}_{\infty,\beta}$ is convex, by general
results concerning the problems of the optimal recovery \cite{Sm}, there
exists a linear optimal method of recovery, that is, there exist
numbers $\alpha_{0},\alpha_{1},\beta_{1},\ldots,
\alpha_{n-1},\beta_{n-1}$ such that
\begin{equation}
\begin{array}{ll}
&\dsup_{f\in \widetilde{K}^{G,\varphi}_{\infty,\beta}}
\|f(0)-\alpha_{0}a_{0}(f)-\sum_{j=1}^{n-1}(\alpha_{j}a_{j}(f)
+\beta_{j}b_{j}(f))\|_{\infty}\\
=&\dsup_{f\in\widetilde{K}^{G,\varphi}_{\infty,\beta},
I_{2n-1}(f)=0}\|f\|_{\infty}=\|\Phi_{n,\beta}^{G,\varphi}\|_{\infty}.
\end{array}\label{dt4-20}
\end{equation}
Let $g$ be an arbitrary function in
$\widetilde{K}^{G,\varphi}_{\infty,\beta}$. For $t\in [0,2\pi)$,
we set $f_{t}(\tau)=g(t+\tau)$. We have $f_{t}\in
\widetilde{K}^{G,\varphi}_{\infty,\beta}$,
$a_{0}(f_{t})=a_{0}(g)$, and
\begin{equation}
\begin{array}{ll}
&a_{j}(f_{t})=a_{j}(g)\cos(jt)+b_{j}(g)\sin(jt),\\[2mm]
&b_{j}(f_{t})=-a_{j}(g)\sin(jt)+b_{j}(g)\cos(jt),
\end{array}
 j=1,2,\ldots;\label{dt4-21}
\end{equation}
therefore setting
$$
{\mathcal L}(I_{2n-1}(g)):=\alpha_{0}a_{0}(g)
+\dsum_{j=1}^{n-1}((\alpha_{j}\cos(jt)-\beta_{j}\sin(jt))a_{j}(g)
+(\alpha_{j}\sin(jt)+\beta_{j}\cos(jt))b_{j}(g)),
$$
from (\ref{dt4-20}) with $\tau=0$ we obtain
$$
|g(t)-{\mathcal L}(I_{2n-1}(g))| \leq
\|\Phi_{n,\beta}^{G,\varphi}\|_{\infty}.
$$
Since $t\in [0,2\pi)$ can be chosen arbitrarily, it follows that
$$
\|g-{\mathcal L}(I_{2n-1}(g))\|_{\infty} \leq
\|\Phi_{n,\beta}^{G,\varphi}\|_{\infty}.
$$
For $g=\Phi_{n,\beta}^{G,\varphi}$ the last inequality turns into
equality, which proves Lemma \ref{Lemma dt4-4}.
\end{proof}

 Let $G$ satisfy Property $B$ and $\widetilde{G}$ be
defined by (\ref{dt4-G1}). Consider a subset of
$\widetilde{K}^{G,\varphi}_{\infty,\beta}$ defined by
\begin{equation}
\widetilde{K}^{G,\varphi}_{\infty,\beta}\bigcap
\mathcal{T}_{n}^{\bot}:
=\left\{f\in\widetilde{K}^{G,\varphi}_{\infty,\beta}:\,
I_{2n-1}(f)=0 \right\}, \label{dt4-g3}
\end{equation}
and let $F$ be the periodic integral of
$f\in\widetilde{K}^{G,\varphi}_{\infty,\beta}\bigcap
\mathcal{T}_{n}^{\bot}$ and satisfy the condition
$\int_{0}^{2\pi}F(t)dt=0$. Then we have $I_{2n-1}(F)=0$. By
(\ref{dt4-19}) we know that
\begin{equation}
\|F\|_\infty\leq\|\Phi_{n,\beta}^{\widetilde{G},\varphi}\|_{\infty},
\label{dt4-gg}
\end{equation}
which together with  Theorem \ref{Theorem dt4-3} gives  an
inequality of Taikov type as follows.

\begin{theorem}[Inequalities of Taikov type ]\label{Theorem dt4-4} Let $n \in{\mathbb
N}$, $G$ satisfy Property $B$ and
$f\in\widetilde{K}^{G,\varphi}_{\infty,\beta}\bigcap
\mathcal{T}_{n}^{\bot}$. Then
$$
\|f\|_q\leq\|\Phi_{n,\beta}^{G,\varphi}\|_q,\quad 1\leq q<\infty.
$$
\end{theorem}

 Since $\Phi_{n,\beta}^{G,\varphi}$ is
$2\pi/n$-periodic, $\Phi_{n,\beta}^{G,\varphi}\in
\widetilde{K}^{G,\varphi}_{\infty,\beta} \bigcap
\mathcal{T}_{n}^{\bot}$. From Theorem \ref{Theorem dt4-4} and
(\ref{dt4-gg}), we obtain

\begin{corollary}\label{Corollary dt4-5} Let $n\in{\mathbb N}$ and $G$ satisfy
Property $B$. Then
$$
\sup_{f\in\widetilde{K}^{G,\varphi}_{\infty,\beta} \bigcap
\mathcal{T}_{n}^{\bot}}\|f\|_q=\|\Phi_{n,\beta}^{G,\varphi}\|_q,
\quad 1\leq q<\infty.
$$
\end{corollary}

\begin{remark}\label{remark dt4-2} Among others, Taikov \cite{Ta} proved that for
all $n\in {\mathbb N}$ and $r=1,2,\ldots$,
\begin{equation}
\sup_{f\in\widetilde{W}_\infty^r \bigcap
\mathcal{T}_{n}^{\bot}}\|f\|_q\leq\|D_{r}\ast h_{n}\|_q, \quad
1\leq q<\infty.\label{dt4-ss}
\end{equation}
The inequality (\ref{dt4-ss}) is now referred to as Taikov's
inequality. For the case  $q=1$, (\ref{dt4-ss}) was also proved by
Turovets in another way. For more details see \cite[p. 172 and p.
207]{Ko}. We \cite{FGSLXH} proved that for all $n\in {\mathbb N}$
and $r=0,1,2,\ldots$,
\begin{equation}
\sup_{f\in\widetilde{H}_{\infty,\beta}^{r}\bigcap
\mathcal{T}_{n}^{\bot}}\|f\|_q=\|D_{r}\ast\varphi_{0}(K_{\beta}\ast
h_{n})\|_q, \quad 1\leq q<\infty. \label{dt4-ss3}
\end{equation}
\end{remark}

\begin{remark}\label{remark dt4-2-0} We conjecture that the
inequality of Taikov type also holds for $G$ being a $NCVD$
kernel, the maim difficulty lies proving the corollary \label{Corollary dt4-2}
is also true for $G$ being a $NCVD$  kernel.
\end{remark}

\section{ Exact Values of n--Widths of
$\widetilde{K}^{G,\varphi}_{\infty,\beta}$}\label{sec4}

\hspace{1ex} In this section, we will solve a uniform minimum norm question
of the classical polynomial perfect splines. Using the result,
we determine the exact values of the Kolmogorov, Gel$'$fand,
linear and information $n$--widths of
$\widetilde{K}^{G,\varphi}_{\infty,\beta}$ in $L_\infty$ for
$\varphi$ being such that
$\widetilde{K}_{\infty,\beta}^{G,\varphi}$ is convex, which have
been obtained by Osipenko \cite{Oc97b} in a different method. And
then we solve a minimum norm question on analogue of the classical
polynomial perfect splines in the space of $L_q$, $1\leq q<\infty$,
which together with some results of
Section 3 gives the exact values of the Gel$'$fand $n$--width, and
the lower estimate of the Kolmogorov $2n$--width of
$\widetilde{K}^{G,\varphi}_{\infty,\beta}$ in $L_{q}$ for $G$
satisfying Property $B$ and $1\leq q<\infty$. \\

 Let $n\in {\mathbb N}$ and let
$$
\overline{\Lambda}_{2n}^{\varphi,\Theta}:
=\{\xi\in\overline{\Lambda}_{2n}:\, \varphi(K_{\beta}\ast
h_{\xi})\perp \Theta\},
$$
\begin{equation}
\widetilde{K}^{G}_{\beta}(\overline{\Lambda}_{2n}^{\varphi,\Theta}):=\{f:\quad
f=a+G\ast\varphi(K_{\beta}\ast
h_{\xi}),\,a\in\Theta,\,\xi\in\overline{\Lambda}_{2n}^{\varphi,\Theta}
\},\label{dt4-14}
\end{equation}
where $G$, $\varphi$, $K_{\beta}$ and $\Theta$ are the same as
those in (\ref{dt4-c}). We will regard $\varphi(K_{\beta}\ast u)$
as $u$ if $\beta=0$.

\begin{lemma}\label{Lemma dt4-5}\ \ For $n\in {\mathbb N}$, we have
\begin{equation}
 \dinf_{f\in\widetilde{K}^{G}_{\beta}(\overline{\Lambda}_{2n}^{\varphi,\Theta})}\|f\|_{\infty}
 =\|G\ast \varphi(K_{\beta}\ast
h_{n})\|_{\infty}=\|\Phi_{n,\beta}^{G,\varphi}\|_{_{\infty}}.
 \label{dt4-15}
\end{equation}
\end{lemma}
\begin{proof} Suppose that there exists
an $a\in\Theta$ and a
$\xi\in\overline{\Lambda}_{2n}^{\varphi,\Theta}$ for which
$$
\|a+G\ast\varphi(K_{\beta}\ast h_{\xi})\|_{\infty}
 <\|G\ast \varphi(K_{\beta}\ast h_{n})\|_{\infty}.
$$
By virtue of (\ref{dt4-11}), we have
\begin{equation}
2n\leq S_{c}((G\ast \varphi(K_{\beta}\ast h_{n}))(\cdot+\alpha)\pm
a\pm(G\ast\varphi(K_{\beta}\ast h_{\xi}))(\cdot)) \label{dt4-16}
\end{equation}
for every $\alpha\in {\mathbb R}$. It follows from the equality
$\varphi(K_{\beta}\ast h_{n})(x+\pi/n)=-\varphi(K_{\beta}\ast
h_{n})(x)$ that $\varphi(K_{\beta}\ast h_{n})$ is
$2\pi/n$--periodic and $\varphi(K_{\beta}\ast h_{n})\perp 1$. From
the fact that $G$ satisfies Property $B$ or is a $NCVD$ kernel,
$a\in\Theta$, $\varphi$ is an odd, continuous and strictly
increasing function, $K_{\beta}$ is a $NCVD$ kernel and
$\xi\in\overline{\Lambda}_{2n}^{\varphi,\Theta}$, it follows that
\begin{equation}
\begin{array}{lll}
&S_{c}((G\ast \varphi(K_{\beta}\ast h_{n}))(\cdot+\alpha)\pm a\pm
(G\ast\varphi(K_{\beta}\ast h_{\xi}))(\cdot))\\
&\leq S_{c}(\varphi((K_{\beta}\ast h_{n})(\cdot+\alpha))\pm
\varphi((K_{\beta}\ast h_{\xi})(\cdot)))\\
&=S_{c}((K_{\beta}\ast h_{n})(\cdot+\alpha)\pm (K_{\beta}\ast
h_{\xi})(\cdot))\\
&\leq S_{c}(h_{n}(\cdot+\alpha)\pm h_{\xi}(\cdot)).
\end{array}
\label{dt4-17}
\end{equation}
Using \cite[Chapt. V, Lemma 4.1]{P}, we obtain the existence of an
$\alpha\in {\mathbb R}$ and $\epsilon=1$ or $-1$ for which
$$
 S_{c}(h_{n}(\cdot+\alpha)-\epsilon h_{\xi}(\cdot))\leq 2(n-1).
$$
So from (\ref{dt4-16}) and (\ref{dt4-17}), we have
$$
2n\leq S_{c}((G\ast \varphi(K_{\beta}\ast
h_{n}))(\cdot+\alpha)-\epsilon a-\epsilon
(G\ast\varphi(K_{\beta}\ast h_{\xi}))(\cdot))\leq 2(n-1)
$$
for some $\alpha\in {\mathbb R}$ and $\epsilon=1$ or $-1$. This is
a contradiction. Lemma \ref{Lemma dt4-5} is proved.
\end{proof}

\begin{theorem}\label{Theorem dt4-5} Let $n\in {\mathbb N}$ and $\varphi$
be such that $\widetilde{K}_{\infty,\beta}^{G,\varphi}$ is convex.
Then
\begin{equation}
\begin{array}{ll}
d_{2n}(\widetilde{K}^{G,\varphi}_{\infty,\beta},L_{\infty})
&=\lambda_{2n}(\widetilde{K}^{G,\varphi}_{\infty,\beta},L_{\infty})
=d^{2n}(\widetilde{K}^{G,\varphi}_{\infty,\beta},L_{\infty})
=i_{2n}(\widetilde{K}^{G,\varphi}_{\infty,\beta},L_{\infty})\\
&=d_{2n-1}(\widetilde{K}^{G,\varphi}_{\infty,\beta},L_{\infty})
=\lambda_{2n-1}(\widetilde{K}^{G,\varphi}_{\infty,\beta},L_{\infty})\\
&=d^{2n-1}(\widetilde{K}^{G,\varphi}_{\infty,\beta},L_{\infty})
=i_{2n-1}(\widetilde{K}^{G,\varphi}_{\infty,\beta},L_{\infty})
=\|\Phi_{n,\beta}^{G,\varphi}\|_{\infty}.
\end{array}
\label{dt4-22}
\end{equation}
\end{theorem}

\begin{proof} First, we will prove the
lower bound for the Kolmogorov $2n$--width. Set
\begin{equation}
\begin{array}{ll}
S^{2n}:&=\{x=(x_1,\ldots,x_{2n+1})\in {\mathbb R}^{2n+1} :\,
\sum_{k=1}^{2n+1}|x_{k}|=2\pi\}, \\
\tau_{0}(x):&=0,\quad \tau_{j}(x):=\sum_{k=1}^{j}|x_{k}|, \quad
j=1,\ldots,2n+1.
\end{array}
\label{dt6-S}
\end{equation}
For each $x\in S^{2n}$, put
\begin{equation}
\begin{array}{ll}
g_{x}(y):&={\rm sgn}\, x_{j},\quad \tau_{j-1}(x)\leq
y<\tau_{j}(x),\quad j=1,\ldots,2n+1 , \\
f_{x}(y):&=(G\ast \varphi(K_{\beta}\ast g_{x}))(y).
\end{array}
\label{dt6-g-f}
\end{equation}

Let $X_{2n}$ be any $2n$--dimensional subspace of $L_q$ ,
$1<q<\infty$. Suppose that $X_{2n}={\rm
span}\{f_{1},\ldots,f_{2n}\}$ and let
$$
f_{x}^{\circ}:=\sum_{j=1}^{2n}a_{j}(x)f_{j}
$$
be the unique best approximation  element to $f_{x}$ from the
subspace $X_{2n}$. If $\Theta$ is not contained in $X_{2n}$ (this
happens only when $G$ satisfies Property $B$), then
$E(\widetilde{K}^{G,\varphi}_{\infty,\beta},X_{2n})=\infty$. We
shall now assume that $\Theta \subset X_{2n}$. Then $\dim\Theta=1$
if $G$ satisfies Property $B$ or  $\dim\Theta=0$ if $G$ is a
$NCVD$ kernel. Assume that $f_{1}$ is a basis of $\Theta$ when $G$
satisfies Property $B$. The mapping
$$
A_{1}(x):=\left\{
\begin{array}{ll}
&(b(x),a_{2}(x),\ldots,a_{2n}(x)),\quad \mbox{if $G$ satisfies Property $B$},\\
&(a_{1}(x),a_{2}(x),\ldots,a_{2n}(x)),\quad \mbox{if $G$ is a
$NCVD$ kernel},
\end{array}
\right.
$$
where
$$
b(x):=\int_{0}^{2\pi}\varphi((K_{\beta}\ast g_{x})(t))\,dt,
$$
is an odd and continuous map of $S^{2n}$ into ${\mathbb R}^{2n}$.
Hence, by Borsuk's theorem (see \cite{Ko}, p. 91) there exists an
$x_{1}^{\star}\in S^{2n}$ for which $A_{1}(x_{1}^{\star})=0$. Then
$g_{x_{1}^{\star}}\in \{h_{\xi}:\xi\in
\overline{\Lambda}_{2n}^{\varphi,\Theta}\}$ and
$a_{i}(x_{1}^{\star})=0$, $i=2,\ldots,2n$ if $G$ satisfies
Property $B$ or $i=1,\ldots,2n$ if $G$ is a $NCVD$ kernel. As
$\|g_{x_{1}^{\star}}\|_{\infty}\leq 1$, $f_{x_{1}^{\star}}\in
\widetilde{K}^{G,\varphi}_{\infty,\beta}$. Therefore, we have
\begin{equation}
\begin{array}{lll}
E(\widetilde{K}^{G,\varphi}_{\infty,\beta},X_{2n})&\geq
\dinf_{g\in X_{2n}}\|f_{x_{1}^{\star}}-g\|_{q}\geq
\|f_{x_{1}^{\star}}-f_{x_{1}^{\star}}^{\circ}\|_{q}\\
&\geq\dinf_{f\in\widetilde{K}^{G}_{\beta}(\overline{\Lambda}_{2n}^{\varphi,\Theta})}\|f\|_{q},
\quad 1<q<\infty. \label{dt4-26}
\end{array}
\end{equation}
Consequently,
$$
d_{2n}(\widetilde{K}^{G,\varphi}_{\infty,\beta},L_{q})\geq
\dinf_{f\in\widetilde{K}^{G}_{\beta}(\overline{\Lambda}_{2n}^{\varphi,\Theta})}\|f\|_{q},
\quad 1<q<\infty.
$$
Passing to the limit $q\rightarrow\infty$ and using Lemma
\ref{Lemma dt4-5}, we obtain
 $$
 d_{2n}(\widetilde{K}^{G,\varphi}_{\infty,\beta},L_{\infty})\geq
\dinf_{f\in\widetilde{K}^{G}_{\beta}(\overline{\Lambda}_{2n}^{\varphi,\Theta})}\|f\|_{\infty}
=\|\Phi_{n,\beta}^{G,\varphi}\|_{\infty}.
$$

\hspace{1ex} Now let us find a lower bound for the Gel$'$fand
$2n$--width. Suppose that
$$
X^{2n}:=\{f\in L_{\infty}:\,\langle l_{j},\, f\rangle=0,\,l_{j}\in
{L_{\infty}}^{*},\,j=1,\ldots,2n\}.
$$
If $\langle l_{j}$, $1\rangle=0,j=1,\ldots,2n$, then
$a\in\widetilde{K}^{G,\varphi}_{\infty,\beta}\cap X^{2n}$ for all
$a\in {\mathbb R}$ (this happens only when $G$ satisfies Property
$B$), which gives
$$
\sup_{f\in\widetilde{K}^{G,\varphi}_{\infty,\beta}\cap
X^{2n}}\|f\|_{\infty}=\infty.
$$
If there exists a $j_{0}\in\{1,\ldots,2n\}$ for which $\langle
l_{j_{0}},1\rangle\neq 0$, then without loss of generality we may
assume that $\langle l_{1},1\rangle\neq 0$. Set
$$
L_{j}:=l_{j}-\frac{\langle l_{j},1\rangle}{\langle
l_{1},1\rangle}l_{1},\quad j=2,\ldots,2n.
$$
For each $x\in S^{2n}$, denote by $A_{2}$ the mapping
$$
A_{2}(x):=\left\{
\begin{array}{ll}
&(b(x),\langle L_{2},f_{x}\rangle,\ldots,\langle
L_{2n},f_{x}\rangle),\quad \mbox{if $G$ satisfies Property $B$},\\
&(\langle l_{1},f_{x}\rangle,,\langle
l_{2},f_{x}\rangle,\ldots,\langle l_{2n},f_{x}\rangle),\quad
\mbox{if $G$ is a $NCVD$ kernel}.
\end{array}
\right.
$$
Since $A_{2}$ is an odd and continuous map of $S^{2n}$ into
${\mathbb R}^{2n}$, by Borsuk's theorem there exists an
$x_{2}^{\star}\in S^{2n}$ for which $A_{2}(x_{2}^{\star})=0$. Then
$g_{x_{2}^{\star}}\in \{h_{\xi}:\xi\in
\overline{\Lambda}_{2n}^{\varphi,\Theta}\}$ and
$$f_{2}^{\star}:=f_{x_{2}^{\star}}-
\frac{\langle l_{1},f_{x_{2}^{\star}}\rangle}{\langle
l_{1},1\rangle}\in X^{2n}\cap
\widetilde{K}^{G,\varphi}_{\infty,\beta}.
$$
Consequently,
$$
\dsup_{f\in\widetilde{K}^{G,\varphi}_{\infty,\beta}\cap
X^{2n}}\|f\|_{\infty} \geq \|f_{2}^{\star}\|_{\infty} \geq
\dinf_{f\in\widetilde{K}^{G}_{\beta}(\overline{\Lambda}_{2n}
^{\varphi,\Theta})}\|f\|_{\infty}
=\|\Phi_{n,\beta}^{G,\varphi}\|_{\infty}.
$$
Therefore,
\begin{equation}
d^{2n}(\widetilde{K}^{G,\varphi}_{\infty,\beta},L_{\infty})\geq
\|\Phi_{n,\beta}^{G,\varphi}\|_{\infty}. \label{dt4-27}
\end{equation}
It follows from \cite[Lemma 1]{Oc97b} and the monotonicity of
$n$-widths that
$$
\begin{array}{ll}
d_{2n}(\widetilde{K}^{G,\varphi}_{\infty,\beta},L_{\infty})&\leq
\lambda_{2n}(\widetilde{K}^{G,\varphi}_{\infty,\beta},L_{\infty})\leq
\lambda_{2n-1}(\widetilde{K}^{G,\varphi}_{\infty,\beta},L_{\infty}),
\\
d_{2n}(\widetilde{K}^{G,\varphi}_{\infty,\beta},L_{\infty})&\leq
d_{2n-1}(\widetilde{K}^{G,\varphi}_{\infty,\beta},L_{\infty})\leq
\lambda_{2n-1}(\widetilde{K}^{G,\varphi}_{\infty,\beta},L_{\infty}),\\
d^{2n}(\widetilde{K}^{G,\varphi}_{\infty,\beta},L_{\infty})&\leq
i_{2n}(\widetilde{K}^{G,\varphi}_{\infty,\beta},L_{\infty})\leq
\lambda_{2n}(\widetilde{K}^{G,\varphi}_{\infty,\beta},L_{\infty})\leq
\lambda_{2n-1}(\widetilde{K}^{G,\varphi}_{\infty,\beta},L_{\infty}),
\\
d^{2n}(\widetilde{K}^{G,\varphi}_{\infty,\beta},L_{\infty})&\leq
d^{2n-1}(\widetilde{K}^{G,\varphi}_{\infty,\beta},L_{\infty})\leq
i_{2n-1}(\widetilde{K}^{G,\varphi}_{\infty,\beta},L_{\infty})\leq
\lambda_{2n-1}(\widetilde{K}^{G,\varphi}_{\infty,\beta},L_{\infty}).
\end{array}
$$
So it remains to find a suitable upper bound for
$\lambda_{2n-1}(\widetilde{K}^{G,\varphi}_{\infty,\beta},L_{\infty})$.
Such an estimate follows from Lemma \ref{Lemma dt4-4}. Theorem
\ref{Theorem dt4-5} is proved.
\end{proof}

\begin{corollary}\label{corollary dt4-6}\ Let $n,\,r\in {\mathbb
N}$. Then
\begin{equation}
\begin{array}{ll}
&d_{2n}(\widetilde{K}^{G,\varphi}_{\infty,\beta},L_{\infty})
=\lambda_{2n}(\widetilde{K}^{G,\varphi}_{\infty,\beta},L_{\infty})
=d^{2n}(\widetilde{K}^{G,\varphi}_{\infty,\beta},L_{\infty})
=i_{2n}(\widetilde{K}^{G,\varphi}_{\infty,\beta},L_{\infty})\\
=&d_{2n-1}(\widetilde{K}^{G,\varphi}_{\infty,\beta},L_{\infty})
=\lambda_{2n-1}(\widetilde{K}^{G,\varphi}_{\infty,\beta},L_{\infty})
=d^{2n-1}(\widetilde{K}^{G,\varphi}_{\infty,\beta},L_{\infty})
=i_{2n-1}(\widetilde{K}^{G,\varphi}_{\infty,\beta},L_{\infty})\\
=&\left\{
\begin{array}{ll}
 &\|D_{r}\ast h_{n}\|_{\infty},\quad \widetilde{K}^{G,\varphi}_{\infty,\beta}
 =\widetilde{W}^{r}_{\infty},\\
 &\|D_{r}\ast (K_{\beta}\ast h_{n})\|_{\infty},\quad \widetilde{K}^{G,\varphi}_{\infty,\beta}
 =\widetilde{h}^{r}_{\infty,\beta},\\
 &\|D_{r}\ast \varphi_{0}(K_{\beta}\ast h_{n})\|_{\infty},\quad
 \widetilde{K}^{G,\varphi}_{\infty,\beta}=\widetilde{H}^{r}_{\infty,\beta}.
\end{array}
\right.
\end{array}
\label{dt4-28}
\end{equation}
\end{corollary}

  Now we proceed to solve a minimum norm question on analogue of the
polynomial perfect splines in the space of $L_q$, $1\leq
q<\infty$. Following  the  method of \cite{Mp,Z,FGSLXH}, we get
the following result which will be used in the lower estimates of
the $n$-widths in Theorem \ref{Theorem dt4-6}, and
\ref{theorem-dt6-2}.

\begin{lemma}\label{Lemma dt4-6}
\ Let $n,\,r\in {\mathbb N}$ and $G$ satisfy Property $B$ . Then
for $1\leq q<\infty$,
\begin{equation}
\begin{array}{ll}
&\dinf_{f\in\widetilde{K}^{G}_{\beta}(\overline{\Lambda}_{2n}
 ^{\varphi,{\mathbb R}})}\|f\|_{q}
=\dmin_{m\in {\mathbb N},m\leq
 n}\|\Phi_{m,\beta}^{G,\varphi}\|_{q}=\|\Phi_{n,\beta}^{G,\varphi}\|_{q}\\
 &=\left\{
\begin{array}{ll}
 &\|D_{r}\ast h_{n}\|_{q},\quad G=D_{r},\,\beta=0,\\
 &\|D_{r}\ast (K_{\beta}\ast h_{n})\|_{q},\quad
 G=D_{r},\,\varphi=\varphi_{1},\,\beta>0,\\
 &\|D_{r}\ast \varphi_{0}(K_{\beta}\ast h_{n})\|_{q},\quad
 G=D_{r},\,\varphi=\varphi_{0},\,\beta>0.
\end{array}
 \right.
 \end{array}
\label{dt4-29}
\end{equation}
\end{lemma}
\begin{proof} For $1\leq q<\infty$, we
follow the  approach  of Zensykbaev \cite{Z},
Micchelli and Pinkus \cite{Mp}, and Pinkus \cite{P79}.
\end{proof}

A compactness argument shows that the infimum in
(\ref{dt4-29}) is attained, i.e., there exists an $a^{\star}\in
{\mathbb R}$ and a $\xi^{\star}\in \overline{\Lambda}_{2n}
 ^{\varphi,{\mathbb R}}$, $\xi^{\star}=(\xi_{1}^{\star},\ldots,\xi_{2m}^{\star}) $,
$m\leq n$, such that
$$
\dinf_{a\in {\mathbb R},\ \xi\in \overline{\Lambda}_{2n}
 ^{\varphi,{\mathbb R}}}\|a+G\ast
\varphi(K_{\beta}\ast h_{\xi})\|^q_{q} =\|a^{\star}+G\ast
\varphi(K_{\beta}\ast h_{\xi^{\star}})\|^q_{q}.
$$

 Using the method of Lagrange$'$s multiplier, we find that
the optimal   $a^{\star}$ and
$\xi^{\star}=(\xi_{1}^{\star},\ldots,\xi_{2m}^{\star})$ must
satisfy the following system of nonlinear equations:
\begin{equation}
\int_{0}^{2\pi}f(x)\,dx=0, \label{hd4}
\end{equation}
\begin{equation}
\begin{array}{lll}
 &\frac{(-1)^{j}}{2\pi^{2}}\dint_{0}^{2\pi}f(x)
 \left[\dint_{0}^{2\pi}G(t){\varphi}'((K_{\beta}\ast
 h_{\xi^{\star}})(x-t))K_{\beta}(x-t-\xi_{j}^{\star})\,dt\right]
dx\\
 +&\frac{\theta(-1)^{j}}{\pi}\dint_{0}^{2\pi}{\varphi}'((K_{\beta}\ast
h_{\xi^{\star}})(t))K_{\beta}(t-\xi_{j}^{\star})\,dt=0,\quad
j=1,\ldots,2m,
\end{array}
 \label{hd5}
\end{equation}
$$
\int_{0}^{2\pi}\varphi((K_{\beta}\ast h_{\xi^{\star}})(t))\,dt=0,
$$
where $\theta$ is the Lagrange$'$s multiplier, and
\begin{equation}
f(x):=q|a^{\star}+G\ast \varphi((K_{\beta}\ast
h_{\xi^{\star}})(x))|^{q-1}{\rm sgn}[a^{\star}+G\ast
\varphi((K_{\beta}\ast h_{\xi^{\star}})(x))].\label{fgs3}
\end{equation}
Since $\varphi$ is a continuous odd and strictly increasing
function, $G$ satisfies Property $B$ and $K_{\beta}$ is $NCVD$, we
get
\begin{equation}
\begin{array}{ll}
&S_{c}(f)= S_c\left(a^{\star}+G\ast \varphi(K_{\beta}\ast
h_{\xi^{\star}})\right) \leq
S_c(\varphi(K_{\beta}\ast h_{\xi^{\star}}))\\
\leq &S_c(K_{\beta}\ast h_{\xi^{\star}})\leq S_c(h_{\xi^*})\leq
2n. \label{fgs2}
\end{array}
\end{equation}
 We first claim that the knots of the  vector
$\xi^\star$ are equidistant which means that $\xi_{j+1}^{\star}
-\xi_{j}^{\star}=\pi/m$, $j=1,\ldots,2m$, i.e.,
$h_{\xi^\star}=h_m$. By translation we may assume that
$\xi^{\star}=(\xi_{1}^{\star},\ldots,\xi_{2m}^{\star})$ satisfies
$0=\xi_{1}^{\star}<\xi_{2}^{\star}<\cdots<\xi_{2m}^{\star}<2\pi$
and
$$
\delta=\xi_{2}^{\star}
=\xi_{2}^{\star}-\xi_{1}^{\star}=\min\{\xi_{i+1}^{\star}
-\xi_{i}^{\star}:i=1,\ldots,2m\},
$$
where $\xi_{2m+1}^{\star}=2\pi$.  Assume that $h_{\xi^{\star}}\neq
h_{m} $. It follows from \cite[Chapt. V, Lemma 4.1]{P} that
$$
S_{c}(h_{\xi^{\star}}(\cdot )+h_{\xi^{\star}}(\cdot+\delta))\leq
2(m-1).
$$
 Since $\varphi$ is a continuous odd and strictly
increasing function, and $K_{\beta}$ is $NCVD$,  we have
\begin{equation}
\begin{array}{lll} & S_{c}(\varphi((K_{\beta}\ast
h_{\xi^{\star}})(\cdot))+\varphi((K_{\beta}\ast
h_{\xi^{\star}})(\cdot+\delta)))\\
=& S_{c}((K_{\beta}\ast h_{\xi^{\star}})(\cdot)+(K_{\beta}\ast
h_{\xi^{\star}})(\cdot+\delta))
\leq S_{c}(h_{\xi^{\star}}(\cdot)+h_{\xi^{\star}}(\cdot+\delta))
\leq 2(m-1).\label{fgs}
\end{array}
\end{equation}
Set
$$
\begin{array}{ll}
&p(x)=a^{\star}+G\ast \varphi((K_{\beta}\ast
h_{\xi^{\star}})(x)),\\
& r(x)=p(x+\delta)=a^{\star}+G\ast \varphi((K_{\beta}\ast
h_{\xi^{\star}})(x+\delta)). \end{array}
$$
\noindent Thus
$$
\begin{array}{ll}
& p(x)+r(x)\\
=&2a^{\star}+\frac{1}{2\pi}
\dint_{0}^{2\pi}G(x-y)[\varphi((K_{\beta}\ast
h_{\xi^{\star}})(y))+ \varphi((K_{\beta}\ast
h_{\xi^{\star}})(y+\delta))]\,dy.
\end{array}$$
From $G$ satisfying Property $B$ and (\ref{fgs}), it follows that
$S_{c}(p+r)\leq 2(m-1)$. Since
$$
{\rm sgn}(a+b)={\rm sgn}\,(|a|^{q-1}{\rm sgn}\,a+|b|^{q-1}{\rm
sgn}\,b)
$$
for every $a$,$b\in {\mathbb R}$ and $q\in (1,\infty)$, it follows
that
$$
S_{c}(|p(\cdot)|^{q-1}{\rm sgn}(p(\cdot))+|r(\cdot)|^{q-1}{\rm
sgn}(r(\cdot)))\leq 2(m-1)
$$
for each $q\in [1,\infty)$.  From (\ref{fgs3}),
$$
\begin{array}{ll}
f(x)&=q|a^{\star}+G\ast \varphi((K_{\beta}\ast
h_{\xi^{\star}})(x))|^{q-1}{\rm sgn}[a^{\star}+G\ast
\varphi((K_{\beta}\ast h_{\xi^{\star}})(x))]\\
&=q|p(x)|^{q-1}{\rm sgn}(p(x)).
\end{array}
$$
Therefore, we have
$$
S_{c}\left(f(\cdot)+f(\cdot+\delta)\right)\leq S_c(p+r)\leq
2(m-1).
$$
Set
$$
P(y)=\frac{1}{2\pi}
\int_{0}^{2\pi}f(x)H(x;y)\,dx+\theta\int_{0}^{2\pi}
{\varphi}'((K_{\beta}\ast h_{\xi^{\star}})(t))K_{\beta}(t-y)\,dt
$$
and
$$
R(y)=\frac{1}{2\pi} \int_{0}^{2\pi}f(x+\delta)H(x;y)\,dx+\theta
\int_{0}^{2\pi}{\varphi}'((K_{\beta}\ast
h_{\xi^{\star}})(t+\delta))K_{\beta}(t-y)\,dt,
$$
where
$$ H(x;y)=\int_{0}^{2\pi}G(t){\varphi}'((K_{\beta}\ast
h_{\xi^{\star}})(x-t))K_{\beta}(x-t-y)\,dt.
$$
 By change of scale and Fubini's theorem, we have
$$
\begin{array}{lll}
&P(y)=\frac{1}{2\pi}\dint_{0}^{2\pi}f(x)\left[\dint_{0}^{2\pi}
G(x-t){\varphi}'((K_{\beta}\ast
h_{\xi^{\star}})(t))K_{\beta}(t-y)\,dt\right]dx\\
+&\theta\dint_{0}^{2\pi}{\varphi}'((K_{\beta}\ast
h_{\xi^{\star}})(t))K_{\beta}(t-y)\,dt\\
=&\frac{1}{2\pi}\dint_{0}^{2\pi}{\varphi}'((K_{\beta}\ast
h_{\xi^{\star}})(t))K_{\beta}(t-y)
\left[\dint_{0}^{2\pi}G(x-t)f(x)\,dx\right]dt\\
+&\theta\dint_{0}^{2\pi}{\varphi}'((K_{\beta}\ast
h_{\xi^{\star}})(t))K_{\beta}(t-y)\,dt\\
=&\frac{1}{2\pi}\dint_{0}^{2\pi}{\varphi}'((K_{\beta}\ast
h_{\xi^{\star}})(t))f_{G}(t)K_{\beta}(t-y)\,dt
+\theta\dint_{0}^{2\pi}{\varphi}'((K_{\beta}\ast
h_{\xi^{\star}})(t))K_{\beta}(t-y)\,dt\\
=&\frac{1}{2\pi}\dint_{0}^{2\pi}{\varphi}'((K_{\beta}\ast
h_{\xi^{\star}})(t))K_{\beta}(t-y)\left(f_{G}(t)+2\pi\theta\right)\,dt\\
=&\frac{1}{2\pi}\dint_{0}^{2\pi}K_{\beta}(y-t){\varphi}'((K_{\beta}\ast
h_{\xi^{\star}})(t))\left(f_{G}(t)+2\pi\theta\right)\,dt,
\end{array}
$$
where
$$
f_{G}(t)=\int_0^{2\pi}G(x-t)f(x)\,dx.
$$
Since ${\varphi}'\geq 0$,  ${\varphi}'$ is continuous on $[-1,1]$,
$G$ satisfies Property $B$ and $K_{\beta}$ is $NCVD$, by
(\ref{fgs2}), we conclude that
\begin{equation}
\begin{array}{ll}
&\widetilde{Z}_{c}(P(\cdot)+R(\cdot))\leq
S_c\left({\varphi}'((K_{\beta}\ast
h_{\xi^{\star}})(\cdot))(f_{G}(\cdot)+f_{G}(\cdot+\delta)+4\pi\theta)\right)\\
\leq & S_c(f_{G}(\cdot)+f_{G}(\cdot+\delta)+4\pi\theta)\leq
S_c(f(\cdot)+f(\cdot+\delta))\leq 2(m-1).\label{G3}
\end{array}
\end{equation}
A simple change of variable argument shows that
$R(y)=P(y+\delta)$, from which we obtain
$$
\widetilde{Z}_{c}(P(\cdot)+P(\cdot+\delta))\leq 2(m-1).
$$
From (\ref{hd5}), we have
$$
P(\xi_{i}^{\star})=0,\quad i=1,\ldots,2m.
$$
By our choice of $\delta$,
$\xi_{i}^{\star}<\xi_{i}^{\star}+\delta\leq \xi_{i+1}^{\star}$,
$i=1,\ldots,2m$, and therefore
$$
S_{c}^{+}\left(P(\xi_{1}^{\star}+\delta),
\ldots,P(\xi_{2n}^{\star}+\delta)\right)=2m .
$$
Thus
$$
S_{c}^{+}\left(P(\xi_{1}^{\star})
+P(\xi_{1}^{\star}+\delta),\ldots,P(\xi_{2n}^{\star})
+P(\xi_{2n}^{\star}+\delta)\right)=2m ,
$$
which implies that
$$
\widetilde{Z}_{c}\left(P(\cdot)+P(\cdot+\delta)\right)\geq 2m.
$$
This is a contradiction, and therefore
$h_{\xi^{\star}}(\cdot)=-h_{\xi^{\star}}(\cdot+\delta)$, i.e.,
$h_{\xi^{\star}}=h_{m}$.

 Now we proceed to show that $a^{\star}=0$. Let
$$
f(x;a):=q|a+G\ast \varphi((K_{\beta}\ast h_{m})(x))|^{q-1}{\rm
sgn}[a+G\ast \varphi((K_{\beta}\ast h_{m})(x))].
$$
Since the constant term is a free variable,
$f(\cdot;a^{\star})\perp 1$. Because  $G\ast
\varphi((K_{\beta}\ast h_{m})(x+\frac{\pi}{m}))=-G\ast
\varphi((K_{\beta}\ast h_{m})(x))$, it follows that
$f(\cdot;0)\perp 1$. Thus $a^{\star}=0$. Since $ h_{m}\in
\{h_{\xi}:\xi\in \overline{\Lambda}_{2n}^{\varphi,{\mathbb R}}\}$
for all $m\in {\mathbb N}$ satisfying $m\leq n$, by Theorem
\ref{Theorem dt4-4} and (\ref{dt4-9}), we get
(\ref{dt4-29}).\hfill\qed

\begin{theorem}\label{Theorem dt4-6}\ Let $n,\,r\in {\mathbb N}$ and $G$
satisfy Property $B$. Then for $1\leq q<\infty$,
\begin{equation}
\begin{array}{ll}
d^{2n}(\widetilde{K}^{G,\varphi}_{\infty,\beta},L_{q})
&=d^{2n-1}(\widetilde{K}^{G,\varphi}_{\infty,\beta},L_{q})
=\|\Phi_{n,\beta}^{G,\varphi}\|_{q}\\
&=\left\{
\begin{array}{ll}
 &\|D_{r}\ast h_{n}\|_{q},\quad G=D_{r},\,\beta=0,\\
 &\|D_{r}\ast (K_{\beta}\ast h_{n})\|_{q},\quad G=D_{r},\,\varphi=\varphi_{1},\,\beta>0,\\
 &\|D_{r}\ast \varphi_{0}(K_{\beta}\ast h_{n})\|_{q},\quad
 G=D_{r},\,\varphi=\varphi_{0},\,\beta>0.
 \end{array}
\right.
\end{array}
\label{dt4-31}
\end{equation}
\end{theorem}
\begin{proof} We begin with the lower
estimate. Assume that $S^{2n}$, $g_{x}$ and $f_{x}$ are the same
as those in the proof of Theorem \ref{Theorem dt4-5} (see
(\ref{dt6-S}) and (\ref{dt6-g-f})). Suppose that $1\leq q<\infty$,
and
$$
X^{2n}:=\{f\in L_{q}:\,\langle l_{j},\, f\rangle=0,\,l_{j}\in
{L_{q}}^*,\,j=1,\ldots,2n\}.
$$
If $\langle l_{j}$, $1\rangle=0,j=1,\ldots,2n$, then
$$
\sup_{f\in \widetilde{K}^{G,\varphi}_{\infty,\beta}\cap
X^{2n}}\|f\|_{q}=\infty.
$$
Therefore, we only need to consider the subspace $X^{2n}$ such
that there exists a $j_{0}\in\{1,\ldots,2n\}$ for which $\langle
l_{j_{0}},1\rangle\neq 0$. Then without loss of generality we may
assume that $\langle l_{1},1\rangle\neq 0$. Set
$$
L_{j}:=l_{j}-\frac{\langle l_{j},1\rangle}{\langle
l_{1},1\rangle}l_{1},\quad j=2,\ldots,2n.
$$
For each $x\in S^{2n}$, denote by $A_{3}$ the mapping
$$
A_{3}(x):=(b(x),\langle L_{2},f_{x}\rangle,\ldots,\langle
L_{2n},f_{x}\rangle),
$$
where
$$
b(x):=\int_{0}^{2\pi}\varphi((K_{\beta}\ast g_{x})(t))dt.
$$
Since $A_{3}$ is an odd and continuous map of $S^{2n}$ into
${\mathbb R}^{2n}$, by Borsuk's theorem, there exists an
$x_{3}^{\star}\in S^{2n}$ for which $A_{3}(x_{3}^{\star})=0$. Then
$g_{x_{3}^{\star}}\in \{h_{\xi}:\xi\in \overline{\Lambda}_{2n}
 ^{\varphi,{\mathbb R}}\}$, and $\langle
L_{i},f_{x_{3}^{\star}}\rangle=0$, $i=2,\ldots,2n$. Thus
$$f_{3}^{\star}:=f_{x_{3}^{\star}}-
\frac{\langle l_{1},f_{x_{3}^{\star}}\rangle}{\langle
l_{1},1\rangle}\in X^{2n}\cap
\widetilde{K}^{G,\varphi}_{\infty,\beta}.
$$
Consequently,
$$
\dsup_{f\in\widetilde{K}^{G,\varphi}_{\infty,\beta}\cap
X^{2n}}\|f\|_{q} \geq \|f_{3}^{\star}\|_{q} \geq
\dinf_{f\in\widetilde{K}^{G}_{\beta}(\overline{\Lambda}_{2n}
 ^{\varphi,{\mathbb R}})}\|f\|_{q}.
$$
Therefore,
\begin{equation}
d^{2n-1}(\widetilde{K}^{G,\varphi}_{\infty,\beta},L_{q})\geq
d^{2n}(\widetilde{K}^{G,\varphi}_{\infty,\beta},L_{q})\geq
\dinf_{f\in\widetilde{K}^{G}_{\beta}(\overline{\Lambda}_{2n}
 ^{\varphi,{\mathbb R}})}\|f\|_{q}.
\label{dt4-h4}
\end{equation}
 Now we turn to the upper estimate. It follows from
Corollary \ref{Corollary dt4-5} that
\begin{equation}
d^{2n}(\widetilde{K}^{G,\varphi}_{\infty,\beta},L_{q})\leq
d^{2n-1}(\widetilde{K}^{G,\varphi}_{\infty,\beta},L_{q})\leq
\|\Phi_{n,\beta}^{G,\varphi}\|_{q}, \label{dt4-h5}
\end{equation}
which together with (\ref{dt4-29}), (\ref{dt4-h4}) and
(\ref{dt4-9}) gives the proof of Theorem \ref{Theorem dt4-6}.
\end{proof}

\begin{remark}\label{remark dt4-4} From \cite[Chapt. V, Theorem 4.11]{P},
it is known that the exact values of the Gel'fand $2n$--width of
 $\widetilde{W}^{r}_{\infty}$ in $L_{q}$, $r\in {\mathbb N}$,
$1\leq q\leq\infty$ is determined, moreover, Sun \cite{S} determined the
exact values of the Gel'fand $n$--widths
$$
d^{2n}(\widetilde{h}_{\infty,\beta},L_{q})
=d^{2n-1}(\widetilde{h}_{\infty,\beta},L_{q})
=\|\phi_{n}^{\beta}\|_{q}, \quad 1\leq q\leq\infty,\ n=1,2,\ldots,
$$
where
\begin{equation}
\phi_{n}^\beta(x)=\frac
4\pi\sum_{\nu=0}^\infty\frac{\cos(2\nu+1)nx}
{(2\nu+1)\cosh(2\nu+1)n\beta}\,. \label{dt4-fgs6}
\end{equation}
Osipenko \cite{Oc97a} obtained the exact values of the Gel'fand
$n$--width of $\widetilde{h}^{r}_{\infty,\beta}$ in $L_{q}$, $r\in
{\mathbb N}$, $1\leq q\leq\infty$, and we \cite{FGSLXH} determined the
exact values of the Gel'fand $n$--width of
$\widetilde{H}^{r}_{\infty,\beta}$ in $L_{q}$, $r\in {\mathbb N}$,
$1\leq q<\infty$.
\end{remark}
 Now we turn to prove the lower estimate of the
Kolmogorov $2n$--Width of
$\widetilde{K}^{G,\varphi}_{\infty,\beta}$
in $L_{q}$ for $G$ satisfying Property $B$ and $1\leq q<\infty$.

\begin{theorem}\label{theorem-dt6-2}\ Let $n\in {\mathbb N}$ and $G$
satisfy Property $B$. Then for $1\leq q<\infty$,
\begin{equation}
d_{2n}(\widetilde{K}^{G,\varphi}_{\infty,\beta},L_{q})
\geq\|\Phi_{n,\beta}^{G,\varphi}\|_{q}. \label{dt6-17}
\end{equation}
\end{theorem}
\begin{proof} We will prove the lower
bound for the Kolmogorov $2n$--width. Let $n\in {\mathbb N}$ and
$1<q<\infty$. Assume that $S^{2n}$, $g_{x}$ and $f_{x}$  are the
same as those in the proof of Theorem \ref{Theorem dt4-5} (see
(\ref{dt6-S})--(\ref{dt6-g-f})). Since the  class of functions
$\widetilde{K}^{G,\varphi}_{\infty,\beta}$ contains all constants,
in order to establish the lower bound of
$d_{2n}(\widetilde{K}^{G,\varphi}_{\infty,\beta},L_{q})$ we only
have to consider the subspace $X_{2n}\subset L_q$ which also
contains the constants. Let $X_{2n}$ be any $2n$-dimensional
subspace of $L_{q}$, $1<q<\infty$, such that $1\in X_{2n}$.
Suppose that $X_{2n}={\rm span}\{f_{1},\ldots,f_{2n}\}$ and
$f_{1}(t)\equiv 1$. Let $a_{1}(x),\ldots,a_{2n}(x)$ be the
coefficients of $f_{1},\ldots,f_{2n}$, respectively, in the unique
best approximation to $f_{x}$ from $X_{2n}$. The mapping
$$
A_{4}(x):=(b(x),a_{2}(x),\ldots,a_{2n}(x)),
$$
where
$$
b(x):=\int_{0}^{2\pi}\varphi((K_{\beta}\ast g_{x})(t))\,dt,
$$
is an odd and continuous map of $S^{2n}$ into ${\mathbb R}^{2n}$.
By Borsuk's theorem there exists an $x_{4}^{\star}\in S^{2n}$ for
which $A_{4}(x_{4}^{\star})=0$. Then $g_{x_{4}^{\star}}\in
\{h_{\xi}:\xi\in \overline{\Lambda}_{2n}^{\varphi,{\mathbb R}}\}$
and $a_{i}(x_{4}^{\star})=0$, $i=2,\ldots,2n$. Then
$f_{x_{4}^{\star}}-a_{1}(x_{4}^{\star})\in
\widetilde{K}^{G}_{\beta}(\overline{\Lambda}_{2n}^{\varphi,{\mathbb
R}})\subset \widetilde{K}^{G,\varphi}_{\infty,\beta}$. Therefore,
if $1<q<\infty$, by Lemma \ref{Lemma dt4-6} , we have
\begin{equation}
\begin{array}{lll}
&d_{2n}(\widetilde{K}^{G,\varphi}_{\infty,\beta},L_{q})
\geq\dsup_{f\in\widetilde{K}^{G,\varphi}_{\infty,\beta}}\inf_{g\in
X_{2n}}\|f-g\|_{q}\geq
\|f_{x_{4}^{\star}}-a_{1}(x_{4}^{\star})\|_{q}\\
&\geq\dinf_{f\in
\widetilde{K}^{G}_{\beta}(\overline{\Lambda}_{2n}^{\varphi,{\mathbb
R}})}\|f\|_{q}=\|\Phi_{n,\beta}^{G,\varphi}\|_{q}, \quad
1<q<\infty. \label{dt6-K1}
\end{array}
\end{equation}
By passing to the limit $q\rightarrow 1$, we obtain the lower
estimate of the Kolmogorov $2n$--width
$d_{2n}(\widetilde{K}^{G,\varphi}_{\infty,\beta},L_{q})$, $1\leq
q<\infty$.  Theorem \ref{theorem-dt6-2} is proved.
\end{proof}

\begin{remark}\label{remark-dt6-2-0}\ By
\cite{FGSLXH} \cite{Ko} \cite{Oc95b}, and\cite{P}, if  $n,\,r\in {\mathbb
N}$, then for
$\widetilde{K}^{G,\varphi}_{\infty,\beta}=\widetilde{W}^{r}_{\infty}$,
$\widetilde{h}^{r}_{\infty,\beta}$ or
$\widetilde{H}^{r}_{\infty,\beta}$,
\begin{equation}
d_{2n}(\widetilde{K}^{G,\varphi}_{\infty,\beta},L_{q})
=\lambda_{2n}(\widetilde{K}^{G,\varphi}_{\infty,\beta},L_{q})
=i_{2n}(\widetilde{K}^{G,\varphi}_{\infty,\beta},L_{q})
=\|\Phi_{n,\beta}^{G,\varphi}\|_{q}, \label{ff1}
\end{equation}
where
$$
\|\Phi_{n,\beta}^{G,\varphi}\|_{q}\\
=\left\{
\begin{array}{ll}
 &\|D_{r}\ast h_{n}\|_{q},\quad \widetilde{K}^{G,\varphi}_{\infty,\beta}
 =\widetilde{W}^{r}_{\infty},\\
 &\|D_{r}\ast (K_{\beta}\ast h_{n})\|_{q},\quad \widetilde{K}^{G,\varphi}_{\infty,\beta}
 =\widetilde{h}^{r}_{\infty,\beta},\\
 &\|D_{r}\ast \varphi_{0}(K_{\beta}\ast h_{n})\|_{q},\quad
 \widetilde{K}^{G,\varphi}_{\infty,\beta}=\widetilde{H}^{r}_{\infty,\beta},
\end{array}
\right.
$$
and $1\leq q<\infty$. So we conjecture that
$$
d_{2n}(\widetilde{K}^{G,\varphi}_{\infty,\beta},L_{q})
 =\lambda_{2n}(\widetilde{K}^{G,\varphi}_{\infty,\beta},L_{q})
=i_{2n}(\widetilde{K}^{G,\varphi}_{\infty,\beta},L_{q})
 =\|\Phi_{n,\beta}^{G,\varphi}\|_{q}
$$
holds for $n\in {\mathbb N}$, $G$ satisfying Property $B$,
$\varphi$ being such that
$\widetilde{K}_{\infty,\beta}^{G,\varphi}$ is convex and $1\leq
q<\infty$. Moreover, we think that all the results proved in this
paper also hold for $G$ being a $ENCVD$ kernel (see Section 3 of Chapter
III in \cite{P}).
\end{remark}


\begin{thebibliography}{99}

\bibitem{Ac56} N. I. Achieser, \textit{Theory of Approximation},
ngar, New York, 1956.

\bibitem{FGSLXH}  Fang Gensun and  Li Xuehua,
\textsl{ Comparison theorems of Kolmogorov type and exact values
of $n$-widths on Hardy-Sobolev classes}, Math. Comp. 75(2006), 241-258.

\bibitem{FGSLXH2} Fang Gensun and Li Xuehua,
\textsl{Optimal Quadrature Problem on Hardy-Sobolev classes},
 Journal of Complexity, 21 (2005),
722-739.

\bibitem{Fi} S. D. Fisher,
\textit{Envelopes, widths, and Landau problems for analytic
functions}, Constr. Approx. \textbf{5} (1989), 171--187.

\bibitem{Fo} W. Forst, \textit{\"{U}ber die Breite von Klassen holomorpher
periodischer funktionen}, J. Approx. Theory, \textbf{19} (1977),
325--331.

\bibitem{Ko} N. Korneichuk,
\textit{Exact Constants in Approximation Theory}, Cambridge
University Press, Cambridge, 1991.

\bibitem{Lo}
G. G. Lorentz, M. V. Golischek, Y. Makovoz, \textit{Constructive
Approximation}, Springer--Verlag, New York, 1993.

\bibitem{Mp} C. A.
Micchelli, A. Pinkus, \textit{Some problems in the approximation
of functions of two variables and $n$--widths of integral
operators}, J. Approx. Theory \textbf{24} (1978), 51--77.

\bibitem{Oc94} K. Yu.
Osipenko, \textit{On $n$--widths, optimal quadrature  formulas,
and optimal recovery of functions analytic in a strip}, Izv. Ross.
Akad. Nauk, Ser. Mat. \textbf{58} (1994), 55--79; English transl.,
in Russian Acad. Sci. Izv. Math. \textbf{45} (1995), 55--78.

\bibitem{Oc94b} K. Yu.
Osipenko, \textit{Inequalities for derivatives of functions
analytic in a strip}, Mat. Zametki \textbf{56} (1994), 114--122;
English transl., in Math. Notes \textbf{56} (1994), 1069--1074.


\bibitem{Oc95b} K. Yu. Osipenko,
\textit{Exact values of $n$--widths and optimal quadratures on
classes of bounded analytic and harmonic functions}, J. Approx.
Theory \textbf{82} (1995), 156--175.

\bibitem{Oc97a} K. Yu. Osipenko,
\textit{Exact $n$--widths of Hardy--Sobolev classes}, Constr.
Approx. \textbf{13} (1997), 17--27.

\bibitem{Oc97b} K. Yu. Osipenko,
\textit{On the precise values of $n$--widths for classes defined
by cyclic variation diminishing operators}, Sbornik Math.
\textbf{188} (1997), 1371--1383.

\bibitem{Ow} K. Yu. Osipenko and K. Wilderotter,
\textit{Optimal information for approximating periodic analytic
functions}, Math. Comput. \textbf{66} (1997), 1579--1592.

\bibitem{P79} A. Pinkus, \textit{On $n$--widths  of periodic
functions}, J. Analyse Math. \textbf{35} (1979), 209--235.

\bibitem{P} A. Pinkus, \textit{$n$--Widths in Approximation
Theory}, Springer--Verlag, Berlin, 1985.

\bibitem{Sm} S.A.Smolyak, \textit{On the Optimal Recovery of
Functions and Functionals of Them}, Kandidat thesis, Moscow State
University, Moscow, 1965.

\bibitem{S} Sun Yongsheng,
\textit{A remark on  Kolmogorov's comparison theorem}, Chin. Ann.
of  Math. (Ser. B), \textbf{7} (1986), 463--467.

\bibitem{Ta} L. V.  Taikov, \textit{The approximation in the mean
of certain classes of periodic functions}, Trudy  Math. Inst.
Stekhlov, \textbf{88} (1967), 61--70.

\bibitem{T} V. M. Tikhomirov, \textit{Diameters
of sets in function spaces and the theory of best approximations},
Uspekhi Mat. Nauk, \textbf{15} (1960), 81--120; English transl. in
Russian Math. Surveys, \textbf{15} (1960), 75--111.

\bibitem{Z} A. A. Zensykbaev, \textit{On the best quadrature formulas
on the class $W^rL_p$}, Dokl. Akad. Nauk SSSR, \textbf{227}
(1976), 277--279; English transl. in Soviet Math. Dokl,
\textbf{17} (1976), 377--380.

\end{thebibliography}
\end{document}